\documentclass[11pt]{article}
\usepackage{mathrsfs}
\usepackage{amsthm}
\usepackage{amssymb}
\usepackage{amsmath}
\usepackage{graphicx}
\usepackage{color}
\usepackage{amsfonts}
\usepackage{float}
\usepackage{cite}
\usepackage[text={140mm,210mm},left=35mm,vmarginratio=1:1]{geometry}
\newtheorem{theorem}{Theorem}[section]

\newtheorem{lemma}[theorem]{Lemma}

\numberwithin{equation}{section}
\normalsize

\begin{document}
\title{\textbf{Large and Moderate Deviation Principles for the SIR Epidemic in a Random Environment}}

\author{Xiaofeng Xue \thanks{\textbf{E-mail}: xfxue@bjtu.edu.cn \textbf{Address}: School of Science, Beijing Jiaotong University, Beijing 100044, China.}\\ Beijing Jiaotong University\\ Yumeng Shen \thanks{\textbf{E-mail}: 14271040@bjtu.edu.cn \textbf{Address}: School of Science, Beijing Jiaotong University, Beijing 100044, China.}\\ Beijing Jiaotong University}

\date{}
\maketitle

\noindent {\bf Abstract:} In this paper, we are concerned with SIR epidemics in a random environment on complete graphs, where every edges are assigned with i.i.d. weights. Our main results give large and moderate deviation principles of sample paths of this model.

\quad

\noindent {\bf Keywords:} large deviation, moderate deviation, SIR, epidemic, random environment.

\section{Introduction}\label{section one}

In this paper, we are concerned with large and moderate deviation principles of the stochastic SIR (Susceptible-Infected-Removed) epidemic in a random environment on the complete graph. First we introduce some basic definitions and notations. For any $n\geq 1$, we use $C_n$ to denote the complete graph with $n$ vertices. For later use, we identify $C_n$ with $\{1,2,\ldots,n\}$, then $C_m$ is a subset of $C_n$ for any $m<n$. Assuming that $\rho$ is a positive random variable such that $Ee^{\alpha\rho}<+\infty$ for some $\alpha>0$, then for any integers $1\leq i<j$, let $\rho(\{i,j\})$ be an independent copy of $\rho$. We further assume that $\{\rho(\{i,j\}):~i\neq j\}$ are independent. For simplicity, we write $\rho(\{i,j\})$ as $\rho(i,j)$, hence $\rho(i,j)=\rho(j,i)$. Note that $\rho(i,j)$ can be considered as an edge weight on the edge connecting $i$ and $j$.

After the edge weights $\{\rho(i,j):~i\neq j\}$ are given, the stochastic SIR model $\{\eta_t^n\}_{t\geq 0}$ on $C_n$ is a continuous-time Markov process with state space $\{0,1,-1\}^{C_n}$, i.e., at each vertex $i\leq n$, there is a spin $\eta(i)$ taking values in $\{1,0,-1\}$. For any $\eta\in \{0,1,-1\}^{C_n}, i\leq n$ and $l\in \{1,0,-1\}$, let $\eta^{i,l}$ be the configuration in $\{1,0,-1\}^{C_n}$ such that
\[
\eta^{i,l}(j)=
\begin{cases}
\eta(j) & \text{~if~}j\neq i,\\
l & \text{~if~}j=i,
\end{cases}
\]
then, the generator function $\Omega_n$ of $\{\eta_t^n\}_{t\geq 0}$ is given by
\[
\Omega_n f(\eta)=\sum_{i=1}^n\sum_{l\in \{1,0,-1\}}q(\eta,i,l)\big[f(\eta^{i,l})-f(\eta)\big]
\]
for sufficiently smooth $f$ on $\{1,0,-1\}^{C_n}$, where
\[
q(\eta,i,l)=
\begin{cases}
1 & \text{~if~}\eta(i)=1 \text{~and~}l=-1,\\
\frac{\lambda}{n}\sum_{j=1}^n\rho(i,j)1_{\{\eta(j)=1\}} & \text{~if~} \eta(i)=0 \text{~and~}l=1,\\
0 & \text{~else},
\end{cases}
\]
where $\lambda$ is a positive parameter called the infection rate while $1_A$ is the indicator function of the event $A$.

Intuitively, $\{\eta_t^n\}_{t\geq 0}$ describes the spread of an epidemic on $C_n$. Vertices in state $1$ are infected and those in state $0$ are susceptible while those in state $-1$ are removed. An infected vertex becomes removed at rate $1$ while a susceptible vertex $i$ is infected by an infected vertex $j$ at rate proportional to the weight $\rho(i,j)$ on the edge connecting $i$ and $j$. A removed vertex stays in its state forever.

When $\rho\equiv 1$, our processes reduce to the classic SIR model. Readers can see References \cite{Anderson1991} and \cite{Britton2017}  for a survey of this classic case. When $\rho$ satisfies
\[
P(\rho=1)=p=1-P(\rho=0)
\]
for some $p\in (0,1)$, our processes reduce to the SIR model on the Erd\"{o}s-R\'{e}nyi graph $G(n,p)$. For basic properties of $G(n,p)$, see Chapter 2 of \cite{Durrett2007}.

The main results of this paper give large and moderate deviation principles for the above SIR model with random edge weights. For the large deviation part, our first motivation is to extend the result about the classic case given in \cite{Pardoux2017}. Our second motivation is to propose an available approach for the proofs of large deviation principles for some special examples of density-dependent Markov chains introduced in \cite{Kurtz1978}. For mathematical details, see Sections \ref{section two}, \ref{section three} and \ref{section four}. For the moderate deviation part, our result is an analogue of the moderate deviation principle given in \cite{Xue2019} for density-dependent Markov chains. The proof of our result follows a similar strategy with that given in \cite{Xue2019}, except for some details modified according to the assumption of i.i.d. edge weights. For mathematical details, see Sections \ref{section two} and \ref{section five}.

\section{Main results}\label{section two}
In this section we give our main results. For later use,  we introduced some notations, definitions and assumptions. Let $(X, \mathcal{F},P)$ be the probability space under which $\{\rho(i,j):1\leq i<j\}$ are defined. For any $\omega\in X$, let $P_{\lambda,n}^\omega$ be the probability measure of the process $\{\eta_t^n\}_{t\geq 0}$ with infection rate $\lambda$ and given edge weights $\{\rho(i,j,\omega):1\leq i<j\leq n\}$, i.e., $P_{\lambda,n}^\omega$ is the quenched measure of the process. We define $P_{\lambda,n}$ as
\[
P_{\lambda,n}(\cdot)=\int_X \big[P_{\lambda,n}^\omega(\cdot)\big] P(d\omega),
\]
i.e., $P_{\lambda,n}(\cdot)$ is the annealed measure of the process. For any $t\geq 0$, we define
\[
S_t^n=\sum_{i=1}^n1_{\{\eta_t^n(i)=0\}} \text{~and~} I_t^n=\sum_{i=1}^n1_{\{\eta_t^n(i)=1\}},
\]
i.e., $S_t^n$ is the number of susceptible vertices while $I_t^n$ is the number of infected vertices at moment $t$.
For given $T_0>0$, we use $\mathcal{D}\big([0,T_0], \mathbb{R}^2\big)$ to denote the Skorokhod space of c\`{a}dl\`{a}g functions $f:[0, T_0]\rightarrow \mathbb{R}^2$. For later use, for any $f\in \mathcal{D}\big([0,T_0], \mathbb{R}^2\big)$, $0\leq t\leq T_0$ and $x\in \mathbb{R}^2$, we consider $f_t, x$ as column vectors and write $f_t, x$ as
\[
f_t=\big(s_t(f), i_t(f)\big)^T \text{~and~} x=(s_x,i_x)^T,
\]
where $T$ is the transposition operator.

For any $f\in \mathcal{D}\big([0,T_0], \mathbb{R}^2\big)$, we define
\[
\|f\|=\sup_{0\leq t\leq T_0}\big\{|s_t(f)|+|i_t(f)|\big\}.
\]
For later use, we define $\mathcal{B}$ as the subset of $\mathcal{D}\big([0,T_0], \mathbb{R}^2\big)$ of $f$ with the following properties:

1. $i_t(f), s_t(f)\geq 0$ for all $0\leq t\leq T_0$.

2. $s_t(f)$ and $i_t(f)+s_t(f)$ are both decreasing with $t$.

3. If $i_u(f)=0$ for some $u$, then $s_t(f)=s_u(f), i_t(f)=0$ for any $t\geq u$.

\quad

Throughout this paper, we adopt the following assumption.

Assumption \textbf{A}: $\{\eta_0^n(i)\}_{i=1}^n$ are independent and identically distributed such that
\[
 P(\eta_0^n(1)=0)=p_0 \text{~and~} P(\eta_0^n(1)=1)=p_1
\]
for some $p_0, p_1$ not depending on $n$ with $p_0, p_1>0$ and $p_0+p_1<1$.

Now we give our rate functions. For any $f\in \mathcal{D}\big([0,T_0], \mathbb{R}^2\big)$, we define
\begin{align*}
I_{dyn}(f)=&\sup_{g\in C^2\big([0,T_0],\mathbb{R}^2\big)}\Big\{f_{T_0}\cdot g_{T_0}-f_0\cdot g_0-\int_0^{T_0}f_t\cdot g_t^\prime dt\\
&-\int_0^{T_0}\big(e^{g_t\cdot l_1}-1\big)i_t(f)+\lambda(E\rho)\big(e^{g_t\cdot l_2}-1\big)s_t(f)i_t(f)dt\Big\},
\end{align*}
where $l_1=(0,-1)^T$, $l_2=(-1,1)^T$, $g_t^\prime=\big(\frac{d}{dt}s_t(g),\frac{d}{dt}i_t(g)\big)^T$ and $x\cdot y$ is the scalar product of $x, y$, i.e., $x\cdot y=s_xs_y+i_xi_y$. For any $x\in \mathbb{R}^2$, we define
\[
I_{ini}(x)=\sup_{y\in \mathbb{R}^2}\Big\{y\cdot x-\log\big(1-p_0-p_1+e^{s_y}p_0+e^{i_y}p_1\big)\Big\}.
\]
For given $T_0>0$, we use $\vartheta^n$ to denote the path of $\{\big(\frac{S_t^n}{n},\frac{I_t^n}{n}\big)^T\}_{0\leq t\leq T_0}$.
Now we give our main result about the large deviation of $\{\eta_t^n\}_{t\geq 0}$.
\begin{theorem}\label{theorem 2.1 main LDP}
Under Assumption \textbf{A}, for any open set $O\subseteq \mathcal{D}\big([0,T_0], \mathbb{R}^2\big)$,
\begin{equation}\label{equ LDP lower}
\liminf_{n\rightarrow+\infty}\frac{1}{n}\log P_{\lambda, n}\big(\vartheta^n\in O\big)\geq -\inf_{f\in O\bigcap \mathcal{B}}\big(I_{dyn}(f)+I_{ini}(f_0)\big),
\end{equation}
while for any closed set $C\subseteq \mathcal{D}\big([0,T_0], \mathbb{R}^2\big)$,
\begin{equation}\label{equ LDP upper}
\limsup_{n\rightarrow+\infty}\frac{1}{n}\log P_{\lambda, n}\big(\vartheta^n\in C\big)\leq -\inf_{f\in C\bigcap \mathcal{B}}\big(I_{dyn}(f)+I_{ini}(f_0)\big).
\end{equation}
\end{theorem}
If $\rho\equiv 1$, then our model reduces to the classic SIR model, the large deviation principle of which is a special case of the main theorem given in \cite{Pardoux2017}. The classic SIR model is an example of density-dependent Markov chains introduced in \cite{Kurtz1978}. For some integer $n\geq 1$, a density-dependent Markov chain
$\{X_t^n\}_{t\geq 0}$ is with state space $\mathbb{Z}^d$ for some $d\geq 1$ and evolves as
\[
X_t^n\rightarrow X_t^n+l \text{~at rate~} nF_l(\frac{X_t^n}{n})
\]
for any $l\in \mathcal{A}$, where $\mathcal{A}$ is a given subset of $\mathbb{R}^d$ and $\{F_l\}_{l\in \mathcal{A}}$ are smooth functions on $\mathbb{R}^d$. For instance, for the classic SIR model $\{(S_t^n, I_t^n)^T\}_{t\geq 0}$, $d=2$, $\mathcal{A}=\{(0,-1)^T, (-1,1)^T\}$ and
\[
F_{(0,-1)^T}(x)=i_x, F_{(-1,1)^T}(x)=\lambda s_xi_x
\]
for any $x\in \mathbb{R}^2$. Large deviation principles of density-dependent Markov chains are given in References \cite{Agazzi2018, Pardoux2017, Schwartz1995} and so on respectively under different assumptions of $\mathcal{A}$ and $\{F_l\}_{l\in \mathcal{A}}$. We think our strategy of the proof of Theorem \ref{theorem 2.1 main LDP} can also be utilized in giving large deviations of some special cases of density-dependent Markov chains which are not included in those given in the above references. For mathematical details, see the remark given at the end of Section \ref{section four}.

To give more clear expressions of  $I_{dyn}$ and $I_{ini}$ , we have the following theorem.
\begin{theorem}\label{theorem 2.2 Idyn Expression}
For $x\in \mathbb{R}^2$, if $I_{ini}(x)<+\infty$, then $s_x, s_y\geq 0$, $s_x+s_y\leq 1$ and
\[
I_{ini}(x)=s_x\log \frac{s_x}{p_0}+i_x\log \frac{i_x}{p_1}+(1-s_x-i_x)\log\frac{1-s_x-i_x}{1-p_0-p_1},
\]
where $0\log 0$ is defined as $0$.

If $f\in \mathcal{B}$ and $I_{dyn}(f)<+\infty$, then $f$ is absolutely continuous and
\[
I_{dyn}(f)=\int_0^{T_0}L_t(f)+(i_t(f)+s_t(f))^\prime+i_t(f)+s_t^\prime(f)+\lambda(E\rho)i_t(f)s_t(f) dt,
\]
where
\begin{align*}
L_t(f)=&-(i_t(f)+s_t(f))^\prime\log(-(i_t(f)+s_t(f))^\prime)+(i_t(f)+s_t(f))^\prime\log i_t(f)\\
&-s_t^\prime(f)\log(-s^\prime_t(f))+s_t^\prime(f)\log\big(\lambda(E\rho)i_t(f)s_t(f)\big).
\end{align*}

\end{theorem}
Note that it is obviously that $P_{\lambda,n}(\vartheta^n\in \mathcal{B})=1$ according to the definitions of $S_t^n$ and $I_t^n$, hence we only care about $I_{dyn}(f)$ for $f\in \mathcal{B}$.

\quad

To give moderate deviations of our processes, we first state a law of large numbers. Let $\{\widehat{x}_t=(\widehat{s}_t,\widehat{i}_t)^T\in \mathbb{R}^2:~0\leq t\leq T_0\}$ be the unique solution to the ODE
\[
\begin{cases}
&\frac{d}{dt}\widehat{x}_t=l_1H_1(\widehat{x}_t)+l_2H_2(\widehat{x}_t), \\
&\widehat{x}_0=(p_0,p_1)^T,
\end{cases}
\]
where $l_1=(0,-1)^T, l_2=(-1,1)^T$ as we have introduced and $H_1(x)=i_x, H_2(x)=\lambda(E\rho)s_xi_x$ for any $x\in \mathbb{R}^2$. Note that it is easy to check that this ODE satisfies Lipschitz's condition and hence has a unique solution.
The following lemma is an analogue of law of large numbers of density-dependent Markov processes given in \cite{Kurtz1978}.
\begin{lemma}\label{lemma 2.3 LLNofSIR}
For any $\epsilon>0$, $\lim_{n\rightarrow+\infty}P_{\lambda,n}\big(\|\vartheta^n-\widehat{x}\|\geq \epsilon\big)=0$.
\end{lemma}
By Lemma \ref{lemma 2.3 LLNofSIR}, $\vartheta^n$ converges to $\widehat{x}$ in probability as $n\rightarrow+\infty$ and hence moderate deviations of our processes are concerned with $\{\frac{(S_t^n,I_t^n)^T-n\widehat{x}_t}{a_n}\}_{0\leq t\leq T_0}$ for any positive sequence $\{a_n\}_{n\geq 1}$ satisfying
$\lim_{n\rightarrow+\infty}\frac{a_n}{n}=0$ and $\lim_{n\rightarrow+\infty}\frac{a_n}{\sqrt{n}}=+\infty$. To give the precise result, we define
\[
b_t=\sum_{i=1}^2l_i(\nabla^TH_i)(\widehat{x}_t)=
\begin{pmatrix}
-\lambda(E\rho)\widehat{i}_t & -\lambda(E\rho)\widehat{s}_t\\
\lambda(E\rho)\widehat{i}_t & \lambda(E\rho)\widehat{s}_t-1
\end{pmatrix}
\]
and
\[
\sigma_t=\sum_{i=1}^2l_iH_i(\widehat{x}_t)l_i^T=
\begin{pmatrix}
\lambda(E\rho)\widehat{i}_t\widehat{s}_t & -\lambda(E\rho)\widehat{i}_t\widehat{s}_t\\
-\lambda(E\rho)\widehat{i}_t\widehat{s}_t & \lambda(E\rho)\widehat{i}_t\widehat{s}_t+\widehat{i}_t
\end{pmatrix},
\]
where $\nabla^T=(\frac{\partial}{\partial s_x},\frac{\partial}{\partial i_x})$. Then we give our rate functions. For any $f\in \mathcal{D}\big([0,T_0], \mathbb{R}^2\big)$, we define
\begin{align*}
J_{dyn}(f)=&\sup_{g\in C^2\big([0,T_0], \mathbb{R}^2\big)}\Big\{f_{T_0}\cdot g_{T_0}-f_0\cdot g_0-\int_0^{T_0}f_t\cdot g_t^\prime dt \\
&-\int_0^{T_0}(b_tf_t)\cdot g_tdt-\frac{1}{2}\int_0^{T_0}g_t^T\sigma_t g_tdt \Big\}.
\end{align*}
For any $x\in \mathbb{R}^2$, we define
$J_{ini}(x)=\sup_{y\in \mathbb{R}^2}\Big\{y\cdot x-\frac{1}{2}y^TM_{_0}y\Big\}$,
where
\[
M_{_0}=
\begin{pmatrix}
p_0(1-p_0) & -p_0p_1\\
-p_0p_1 & p_1(1-p_1)
\end{pmatrix}.
\]
For given positive sequence $\{a_n\}_{n\geq 1}$ satisfying
$\lim_{n\rightarrow+\infty}\frac{a_n}{n}=0$ and $\lim_{n\rightarrow+\infty}\frac{a_n}{\sqrt{n}}=+\infty$, we denote by $\nu^n$ the path of $\{\frac{(S_t^n, I_t^n)^T-n\widehat{x}_t}{a_n}\}_{0\leq t\leq T_0}$. Now we give our moderate deviations.
\begin{theorem}\label{Theorem 2.4 MDPofSIR}
Under Assumption \textbf{A}, for any open set $O\subseteq \mathcal{D}\big([0,T_0], \mathbb{R}^2\big)$,
\begin{equation}\label{equ MDP lower bound}
\liminf_{n\rightarrow+\infty}\frac{n}{a_n^2}\log P_{\lambda,n}\big(\nu^n\in O\big)\geq-\inf_{f\in O}(J_{dyn}(f)+J_{ini}(f_0)),
\end{equation}
while for any closed set $C\subseteq \mathcal{D}\big([0,T_0], \mathbb{R}^2\big)$,
\begin{equation}\label{equ MDP upper bound}
\limsup_{n\rightarrow+\infty}\frac{n}{a_n^2}\log P_{\lambda,n}\big(\nu^n\in C\big)\leq-\inf_{f\in C}(J_{dyn}(f)+J_{ini}(f_0)).
\end{equation}
\end{theorem}
Theorem \ref{Theorem 2.4 MDPofSIR} is an analogue of the main result given in \cite{Xue2019}, where moderate deviations of density-dependent Markov chains are investigated. As an application, the moderate deviation of classic SIR model with deterministic initial condition can be given directly according to the main result in \cite{Xue2019}, Theorem \ref{Theorem 2.4 MDPofSIR} is an extension of which to the case where i.i.d weights are assigned on every edges.

To give more clear expressions of $J_{dyn}$ and $J_{ini}$, we have the following theorem.
\begin{theorem}\label{theorem 2.5 JdynEpress}
For any $x\in \mathbb{R}^2$, $J_{ini}(x)=\frac{1}{2}x^TM_{_0}^{-1}x$. For any $f\in\mathcal{D}\big([0,T_0], \mathbb{R}^2\big)$, if $J_{dyn}(f)<+\infty$, then
$f$ is absolutely continuous and
\[
J_{dyn}(f)=\frac{1}{2}\int_0^{T_0}(f_t^\prime-b_tf_t)^T\sigma_t^{-1}(f_t^\prime-b_tf_t) dt.
\]
\end{theorem}
Note that it is easy to check that $M_{_0}$ and $\sigma_t$ are invertible according to their definitions.

The proof of Theorem \ref{theorem 2.1 main LDP} is divided into Sections \ref{section three} and \ref{section four} while an outline of the proof of Theorem \ref{Theorem 2.4 MDPofSIR} is given in Section \ref{section five}. In both proofs, an exponential martingale will be introduced and a generalized version of Girsanov's theorem given in \cite{Schuppen1974} will be utilized. The strategy of our proofs is inspired by those introduced in \cite{Gao2003} and \cite{Xue2019}.

As a preparation for the proof of Equation \eqref{equ LDP lower}, Theorem \ref{theorem 2.2 Idyn Expression} is proved at the beginning of Section \ref{section three}. The core idea of the proof of Theorem \ref{theorem 2.2 Idyn Expression} is to show that $I_{dyn}(f)<+\infty$ implies that there exists $\psi$ such that $f$ is the solution to the ODE $f_t^\prime=l_1e^{\psi_t\cdot l_1}H_1(f_t)+l_2e^{\psi_t\cdot l_2}H_2(f_t)$. The proof of Theorem \ref{theorem 2.5 JdynEpress} is given at the beginning of Section \ref{section five}, where Cauchy-Schwartz's inequality and Riesz's representation theorem will be utilized.

\section{The proof of Equation \eqref{equ LDP lower}}\label{section three}
In this section, we give the proof of Equation \eqref{equ LDP lower}. As a preparation, we first give the proof of Theorem \ref{theorem 2.2 Idyn Expression}. For simplicity, we define
\begin{align*}
\Phi_f(g)=&f_{T_0}\cdot g_{T_0}-f_0\cdot g_0-\int_0^{T_0}f_t\cdot g_t^\prime dt \\
&-\int_0^{T_0}\big(e^{g_t\cdot l_1}-1\big)i_t(f)+\lambda(E\rho)\big(e^{g_t\cdot l_2}-1\big)s_t(f)i_t(f)dt
\end{align*}
for any $g\in C^2\big([0,T_0], \mathbb{R}^2\big)$ and $f\in \mathcal{D}\big([0,T_0], \mathbb{R}^2\big)$, i.e.,
\[
I_{dyn}(f)=\sup_{g\in C^2\big([0,T_0], \mathbb{R}^2\big)}\{\Phi_f(g)\}.
\]

\proof[Proof of Theorem \ref{theorem 2.2 Idyn Expression}]

For $x\in \mathbb{R}^2$ with $I_{ini}(x)<+\infty$, if $s_x<0$, then
\begin{align*}
I_{ini}(x)&\geq \sup_{y:i_y=0}\big\{y\cdot x-\log\big(1-p_0-p_1+e^{s_y}p_0+e^{i_y}p_1\big)\big\}\\
&=\lim_{s_y\rightarrow-\infty}\big[s_xs_y-\log\big(1-p_0+e^{s_y}p_0\big)\big]=+\infty,
\end{align*}
which is contradictory. Hence, $s_x\geq 0$. For the same reason, $s_y\geq 0$. If $s_x+s_y>1$, then
\begin{align*}
I_{ini}(x)&\geq \sup_{y:i_y=s_y}\big\{y\cdot x-\log\big(1-p_0-p_1+e^{s_y}p_0+e^{i_y}p_1\big)\big\}\\
&\geq \lim_{c\rightarrow+\infty} \big[c(s_x+i_x)-\log\big(1+(e^c-1)(p_0+p_1)\big)\big].
\end{align*}
Let $\beta_1(c)=c(s_x+i_x)-\log\big(1+(e^c-1)(p_0+p_1)\big)$, then $\lim_{c\rightarrow+\infty}\frac{d}{dc}\beta_1(c)=s_x+s_y-1>0$, hence
\[
I_{ini}(x)\geq \lim_{c\rightarrow+\infty} \big[c(s_x+i_x)-\log\big(1+(e^c-1)(p_0+p_1)\big)\big]=+\infty,
\]
which is contradictory. Hence, $s_x+s_y\leq 1$. Let $\beta_2(y)=y\cdot x-\log\big(1-p_0-p_1+e^{s_y}p_0+e^{i_y}p_1$. When $s_x>0, s_y>0$ and $s_x+s_y<1$, since
$\beta_2(y)=y\cdot x-\log\big(1-p_0-p_1+e^{s_y}p_0+e^{i_y}p_1\big)$ is concave with respect to each coordinate $s_y$ and $i_y$, $\beta_2$ gets its maximum at $y_0$ given by $\frac{\partial}{\partial s_y}\beta_2(y_0)=\frac{\partial}{\partial i_y}\beta_2(y_0)=0$, i,e,
\[
s_{y_0}=\log \big[\frac{s_x(1-p_0-p_1)}{p_0(1-s_x-s_y)}\big] \text{~and~} i_{y_0}=\log\big[\frac{i_x(1-p_0-p_1)}{p_1(1-s_x-s_y)}\big].
\]
Hence,
\[
I_{ini}(x)=\beta_2(y_0)=s_x\log \frac{s_x}{p_0}+i_x\log \frac{i_x}{p_1}+(1-s_x-i_x)\log\frac{1-s_x-i_x}{1-p_0-p_1}.
\]
The proof of $I_{ini}(x)=s_x\log \frac{s_x}{p_0}+i_x\log \frac{i_x}{p_1}+(1-s_x-i_x)\log\frac{1-s_x-i_x}{1-p_0-p_1}$ for the case where $s_xi_x(1-s_x-i_x)=0$ is similar. We omit the details.

For $f\in \mathcal{B}$ with $I_{dyn}(f)<+\infty$, if $s_t(f)$ is not absolutely continuous, then there exists $\epsilon>0$ such that for any integer $n\geq 1$, there exists $0\leq a_{1,n}<b_{1,n}<a_{2,n}<b_{2,n}<\ldots<a_{{k_n},n}<b_{{k_n},n}\leq T_0$ such that $\sum_{i=1}^{k_n}(b_{i,n}-a_{i,n})\leq \frac{1}{n}$ and
\[
\sum_{i=1}^{k_n}|s_{b_{i,n}}(f)-s_{a_{i,n}}(f)|=-\sum_{i=1}^{k_n}(s_{b_{i,n}}(f)-s_{a_{i,n}}(f))\geq \epsilon.
\]
For any $m>0$, let $K_t^{m,n}$ be defined as $K_s^{m,n}=-m$ when $s\in [a_{i,n},b_{i,n})$ for some $1\leq i\leq k_n$ and $K_s^{m,n}=0$ otherwise. For given $m,n$ and any $l\geq 1$, let $g^l\in C^2\big([0,T_0], \mathbb{R}^2\big)$ such that $i_t(g^l)=0$ and $\lim_{l\rightarrow+\infty}s_t(g^l)=K_t^{m,n}$ for all $0\leq t\leq T_0$. For $f\in \mathcal {B}$, $s(f)$ and $i(f)+s(f)$ are both decreasing and hence are both bounded variation functions. Then,
\[
f_{T_0}\cdot g^l_{T_0}-f_0\cdot g^l_0-\int_0^{T_0}f_t\cdot (g^l_t)^\prime dt=\int_0^{T_0}s_t(g^l)ds_t(f)
\]
and therefore
\begin{align*}
\lim_{l\rightarrow+\infty}\Phi_f(g^l)\geq \int_0^{T_0}K_t^{m,n}ds_t(f)-\frac{1}{n}\lambda(E\rho)e^m\|f\|^2
\geq m\epsilon-\frac{1}{n}\lambda(E\rho)e^m\|f\|^2.
\end{align*}
Then, $I_{dyn}(f)\geq \lim_{l\rightarrow+\infty}\Phi_f(g^l)\geq m\epsilon-\frac{1}{n}\lambda(E\rho)e^m\|f\|^2$.
Since $n$ is arbitrary, $I_{dyn}(f)\geq m\epsilon$. Since $m$ is arbitrary, $I_{dyn}(f)=+\infty$, which is contradictory. Therefore, $s_t(f)$ is absolutely continuous. The proof of absolute continuity of $i_t(f)+s_t(f)$ is similar. Only one detail should be modified that we let $s_t(g^l)=i_t(g^l)$ such that $\lim_{l\rightarrow+\infty}s_t(g^l)=K_t^{m,n}$ in this case. We omit the details. As a result, $s_t(f), i_t(f)+s_t(f)$ are both absolutely continuous and hence $f$ is absolutely continuous.

For $f\in \mathcal{B}$ with $I_{dyn}(f)<+\infty$, since $f$ is absolutely continuous,
\[
\Phi_f(g)=\int_0^{T_0} g_t\cdot f_t^\prime dt-\int_0^{T_0}\big(e^{g_t\cdot l_1}-1\big)i_t(f)+\lambda(E\rho)\big(e^{g_t\cdot l_2}-1\big)s_t(f)i_t(f)dt.
\]
For $f\in \mathcal{B}$, $s_t^\prime(f), \big(i_t(f)+s_t(f)\big)^\prime\leq 0$ while $i_t(f)=0$ implies that $\big(i_t(f)+s_t(f)\big)^\prime=0$ and $s_t(f)=0$ implies that $s_t^\prime(f)=0$, hence there exists $h_1(t), h_2(t)\in [0,+\infty)$ such that they are the solution of
\[
\begin{cases}
s_t^\prime(f)&=-h_2(t)\lambda(E\rho)i_t(f)s_t(f),\\
i_t^\prime(f)&=-h_1(t)i_t(f)+h_2(t)\lambda(E\rho)i_t(f)s_t(f).
\end{cases}
\]
Let $h_1^n\in C^2\big([0,T_0], \mathbb{R}\big)$ such that $h_1^n(t)>0$ for all $t\in [0, T_0]$ while
\[
\lim_{n\rightarrow+\infty}\int_0^{T_0}|h_1^n(t)-h_1(t)|dt=0
\]
and
$h_2^n\in C^2\big([0,T_0], \mathbb{R}\big)$ such that $h_2^n(t)>0$ for all $t\in [0, T_0]$ while
\[
\lim_{n\rightarrow+\infty}\int_0^{T_0}|h_2^n(t)-h_2(t)|dt=0,
\]
then we define $\psi^n\in C^2\big([0,T_0], \mathbb{R}^2\big)$ such that $s_t(\psi^n)$ and $i_t(\psi^n)$ satisfies
$h_1^n(t)=\exp\big\{-i_t(\psi^n)\big\}$ and $h_2^n(t)=\exp\big\{i_t(\psi^n)-s_t(\psi^n)\big\}$ for all $t\in [0, T_0]$.
As a result,
\begin{align*}
I_{dyn}(f) \geq &\lim_{n\rightarrow+\infty}\Phi_f(\psi^n)  \\
=&\lim_{n\rightarrow+\infty}\int_0^{T_0} \psi^n_t\cdot f_t^\prime dt\\
&-\int_0^{T_0}h_1(t)i_t(f)-i_t(f)+\lambda(E\rho)h_2(t)s_t(f)i_t(f)-\lambda(E\rho)s_t(f)i_t(f)dt \\
=& \int_0^{T_0}L_t(f)+(i_t(f)+s_t(f))^\prime+i_t(f)+s_t^\prime(f)+\lambda(E\rho)i_t(f)s_t(f) dt
\end{align*}
according to the definition of $\psi^n$ and $L_t(f)$.

On the other hand, for any $g\in C^2\big([0,T_0], \mathbb{R}\big)$,
\begin{align*}
&\Phi_f(g)
\leq \int_0^{T_0} \sup_{\theta\in \mathbb{R}^2} \big\{\theta\cdot f_t^\prime-\big(e^{\theta\cdot l_1}-1\big)i_t(f)-\lambda(E\rho)\big(e^{\theta\cdot l_2}-1\big)s_t(f)i_t(f)\big\}dt \\
&=\int_0^{T_0}L_t(f)+(i_t(f)+s_t(f))^\prime+i_t(f)+s_t^\prime(f)+\lambda(E\rho)i_t(f)s_t(f) dt
\end{align*}
and hence
\[
I_{dyn}(f)\leq \int_0^{T_0}L_t(f)+(i_t(f)+s_t(f))^\prime+i_t(f)+s_t^\prime(f)+\lambda(E\rho)i_t(f)s_t(f) dt.
\]
Therefore,
\[
I_{dyn}(f)=\int_0^{T_0}L_t(f)+(i_t(f)+s_t(f))^\prime+i_t(f)+s_t^\prime(f)+\lambda(E\rho)i_t(f)s_t(f) dt
\]
and the proof is complete.

\qed

According to a non-rigorous mean-field analysis, $S_t^n\rightarrow S_t^n-1$ at rate
\[
\frac{\lambda}{n}\sum_{i:\eta_t(i)=0}\sum_{j:\eta_t(j)=1}\rho(i,j)\approx \frac{\lambda}{n} S_t^nI_t^n(E\rho).
\]
To give this mean-field analysis a rigorous description, we define
\[
\gamma(C,D)=\sum_{i\in C}\sum_{j\in D}\rho(i,j)
\]
for any $C,D\subseteq C_n$ such that $C\bigcap D=\emptyset$ and
\[
\delta_n=\sup\Big\{\frac{\big|\gamma(C,D)-|C||D|(E\rho)\big|}{n^2}:~C,D\subseteq C_n, C\bigcap D=\emptyset\Big\},
\]
where $|C|$ is the cardinality of $C$. Then, we have the following lemma.
\begin{lemma}\label{lemma 3.1}
For any $\epsilon>0$,
\[
\lim_{n\rightarrow+\infty}\frac{1}{n}\log P(\delta_n>\epsilon)=-\infty.
\]
\end{lemma}

\proof

According to Markov's inequality, for any $\theta>0$ and $C,D\in C_n$ such that $C\bigcap D=\emptyset$,
\[
P\big(\gamma(C,D)-|C||D|(E\rho)\geq n^2\epsilon\big)\leq e^{-\theta n^2\epsilon}\big[Ee^{\theta(\rho-E\rho)}\big]^{|C||D|}.
\]
According to Jensen' inequality, $Ee^{\theta(\rho-E\rho)}\geq e^{\theta E(\rho-E\rho)}=1$ and hence
\[
P\big(\gamma(C,D)-|C||D|(E\rho)\geq n^2\epsilon\big)\leq e^{-\theta n^2\epsilon}\big[Ee^{\theta(\rho-E\rho)}\big]^{n^2}
=\Big[e^{-\theta\epsilon}Ee^{\theta(\rho-E\rho)}\Big]^{n^2}.
\]
According to our assumption of $\rho$, $e^{-\theta\epsilon}Ee^{\theta(\rho-E\rho)}$ is well-defined and differentiable for $\theta\in (-\infty,\alpha)$. Since $e^{-0\epsilon}Ee^{0(\rho-E\rho)}=1$ and
\[
\frac{d}{d\theta}e^{-\theta\epsilon}Ee^{\theta(\rho-E\rho)}\Big|_{\theta=0}=-\epsilon<0,
\]
there exists $\theta_1>0$ such that $e^{-\theta_1\epsilon}Ee^{\theta_1(\rho-E\rho)}<1$ and
\[
P\big(\gamma(C,D)-|C||D|(E\rho)\geq n^2\epsilon\big)\leq \big[e^{-\theta_1\epsilon}Ee^{\theta_1(\rho-E\rho)}\big]^{n^2}.
\]
Follows from a similar analysis, there exists $\theta_2>0$ such that $e^{-\theta_2\epsilon}Ee^{-\theta_2(\rho-E\rho)}<1$ and
\[
P\big(\gamma(C,D)-|C||D|(E\rho)\leq-n^2\epsilon\big)\leq \big[e^{-\theta_2\epsilon}Ee^{-\theta_2(\rho-E\rho)}\big]^{n^2}.
\]
Therefore, there exists $\theta_3>0$ such that
\[
P\big(\big|\gamma(C,D)-|C||D|(E\rho)\big|\geq n^2\epsilon\big)\leq 2e^{-\theta_3n^2}
\]
for any $C,D\subseteq C_n$ such that $C\bigcap D=\emptyset$. Since the number of subsets of $C_n$ is $2^n$,
\[
P\big(\delta_n\geq\epsilon\big)\leq 2e^{-\theta_3n^2}4^n
\]
and hence
\[
\lim_{n\rightarrow+\infty}\frac{1}{n}\log P(\delta_n>\epsilon)\leq \log4-\lim_{n\rightarrow +\infty}n\theta_3=-\infty.
\]

\qed

Our strategy of the proof of Equation \eqref{equ LDP lower} is inspired by those introduced in \cite{Gao2003} and \cite{Xue2019}, where an exponential martingale will be introduced. To give this martingale, we recall some properties of Markov processes. Let $\Omega_n$ be generator of  $\{\eta_t^n\}_{t\geq 0}$ defined as in Section \ref{section one} and $C^{2,1}\big([0, T_0]\times \{1,0,-1\}^{C_n}\big)$ be the set of functions $f:[0, T_0]\times \{1,0,-1\}^{C_n}\rightarrow \mathbb{R}$ such that $f(t,\eta)$ has continuous second-order partial derivative with respect to the coordinate $t$ and has continuous partial derivative with respect to the coordinate $\eta(i)$ for all $1\leq i\leq n$, then $\{\mathcal{M}_t(f)\}_{0\leq t\leq T_0}$ defined as
\begin{equation}\label{equ 3.0}
\mathcal{M}_t(f)=f(t,\eta_t^n)-f(0,\eta_0^n)-\int_0^t (\frac{\partial}{\partial u}+\Omega_n)f(u, \eta^n_u)du
\end{equation}
is a martingale for any $f\in C^{2,1}\big([0, T_0]\times \{1,0,-1\}^{C_n}\big)$ and
\begin{equation}\label{equ 3.1}
<\mathcal{M}(f_1), \mathcal{M}(f_2)>_t=\int_0^{t}\Omega_n(f_1f_2)-f_1\Omega_nf_2-f_2\Omega_n f_1du
\end{equation}
for any $f_1,f_2\in C^{2,1}\big([0, T_0]\times \{1,0,-1\}^{C_n}\big)$.

For any $g\in C^2\big([0,T_0], \mathbb{R}^2\big)$, let
\[
f_g(t,\eta)=i_t(g)\sum_{i=1}^n1_{\{\eta(i)=1\}}+s_t(g)\sum_{i=1}^n1_{\{\eta(i)=0\}},
\]
then $f_g(t,\eta_t^n)=i_t(g)I_t^n+s_t(g)S_t^n=g_t\cdot \big(S_t^n, I_t^n\big)^T$. We further define
\[
H_g(t,\eta_t^n)=e^{f_g(t,\eta_t^n)} \text{~and~}\Lambda_t^n(g)=\frac{H_g(t,\eta_t^n)}{H_g(0,\eta_0^n)}\exp\Big(-\int_0^t\frac{(\frac{\partial}{\partial u}+\Omega_n)H_g(u,\eta_u^n)}{H_g(u,\eta^n_u)}du\Big),
\]
then we have the following lemma.
\begin{lemma}\label{lemma 3.2 LambdaIsMartingale}
For any $\omega\in X$ and $g\in C^2\big([0,T_0], \mathbb{R}^2\big)$, $\{\Lambda_t^n(g)\}_{0\leq t\leq T_0}$ is a martingale with expectation $1$ under the quenched measure
$P_{\lambda,n}^\omega$.

\end{lemma}

\proof

According to Ito's formula,
\begin{equation}\label{equ 3.2}
d\Lambda_{t}^n(g)=\frac{\exp\Big(-\int_0^t\frac{(\frac{\partial}{\partial u}+\Omega_n)H_g(u,\eta_u^n)}{H_g(u,\eta^n_u)}du\Big)}{H_g(0,\eta_0^n)}d\mathcal{M}_{t}(H_g),
\end{equation}
where $\mathcal{M}_t(H_g)$ is defined as in Equation \eqref{equ 3.0} and hence $\{\mathcal{M}_t(H_g)\}_{0\leq t\leq T_0}$ is a martingale.
Therefore, $\{\Lambda_t^n(g)\}_{0\leq t\leq T_0}$ is a local martingale. Since $S_t^n, I_t^n\leq n$ for any $t\geq 0$, $\{\Lambda_t^n(g)\}_{0\leq t\leq T_0}$ are uniformly bounded, which ensures that this local martingale is a martingale.

\qed

By Lemma \ref{lemma 3.2 LambdaIsMartingale}, we define $\widehat{P}_{\lambda,n}^{\omega,g}$ as the quenched measure such that
\[
\frac{d\widehat{P}_{\lambda,n}^{\omega,g}}{dP_{\lambda,n}^\omega}=\Lambda_{T_0}^n(g)
\]
for any $\omega\in X$ and $g\in C^2\big([0,T_0], \mathbb{R}^2\big)$. We further define $\widehat{P}^g_{\lambda,n}$ as the annealed measure such that
\[
\widehat{P}^g_{\lambda,n}(\cdot)=\int_{X}\widehat{P}_{\lambda,n}^{\omega,g}(\cdot)P(d\omega).
\]
For $a\geq 0$, we use $\lfloor a\rfloor$ to denote the largest integer not exceeding $a$. For $x,y\geq 0$ such that $x+y\leq 1$, we define
\[
\widehat{P}^{g,x,y}_{\lambda,n}(\cdot)=\widehat{P}^g_{\lambda,n}\Big(\cdot\Big|S_0^n=\lfloor nx\rfloor, I_0^n=\lfloor ny\rfloor\Big).
\]
Then, we have the following lemma, which is crucial for the proof of Equation \eqref{equ LDP lower}.

\begin{lemma}\label{lemma 3.3 LLNunderChangeMeasure}
For any $x,y\geq 0$ such that $x+y\leq 1$ and any $g\in C^2\big([0,T_0], \mathbb{R}^2\big)$, $\{\big(\frac{S_t^n}{n},\frac{I_t^n}{n}\big)^T\}_{0\leq t\leq T_0}$
converges in $\widehat{P}^{g,x,y}_{\lambda,n}$-probability to the solution
\[\{\widetilde{x}_t=(\widetilde{s}_t,\widetilde{i}_t)^T:~0\leq t\leq T_0\}\]
to the ODE
\[
\begin{cases}
&\frac{d}{dt}\widetilde{s}_t=-e^{i_t(g)-s_t(g)}\lambda(E\rho)\widetilde{s}_t\widetilde{i}_t,\\
&\frac{d}{dt}\widetilde{i}_t=-e^{-i_t(g)}\widetilde{i}_t+e^{i_t(g)-s_t(g)}\lambda(E\rho)\widetilde{s}_t\widetilde{i}_t, \\
&(\widetilde{s}_0,\widetilde{i}_0)=(x,y).
\end{cases}
\]
\end{lemma}

\proof

For $0\leq t\leq T_0$, we define
\[
\widehat{\mathcal{M}}_t(H_g)=\int_0^t \frac{1}{H_g(u-,\eta_{u-}^n)}d\mathcal{M}_{u}(H_g),
\]
then by Equation \eqref{equ 3.2},
\begin{equation}\label{equ 3.4}
d\Lambda_t^n(g)=\Lambda_{t-}^n(g)d\widehat{\mathcal{M}}_t(H_g).
\end{equation}
According to Equation \eqref{equ 3.4} and Theorem 3.2 of \cite{Schuppen1974}, which is a generalized version of Girsanov's thoerem, for any martingale $\{M_t\}_{0\leq t\leq T_0}$ under $P_{\lambda,n}^\omega$,
\[
\{\widetilde{M}_t=M_t-<M,\widehat{\mathcal{M}}(H_g)>_t:~0\leq t\leq T_0\}
\]
is a martingale under $\widehat{P}_{\lambda,n}^{\omega,g}$ and $[\widetilde{M}, \widetilde{M}]=[M, M]$ under both $P_{\lambda,n}^\omega$ and $\widehat{P}_{\lambda,n}^{\omega,g}$.

Let $f_1(\eta)=\sum_{i=1}^n1_{\{\eta(i)=0\}}$ and $f_2(\eta)=\sum_{i=1}^n1_{\{\eta(i)=1\}}$, then, as we have recalled, $\{\mathcal{M}_t(f_1)\}_{0\leq t\leq T_0}$ and $\{\mathcal{M}_t(f_2)\}_{0\leq t\leq T_0}$ are martingales under $P_{\lambda,n}^\omega$, where
\[
\mathcal{M}_t(f_1)=S_t^n-\lfloor nx\rfloor -\int_0^t \Omega_nf_1(\eta_u^n)du
\]
and
\[
\mathcal{M}_t(f_2)=I_t^n-\lfloor ny\rfloor-\int_0^t \Omega_nf_2(\eta_u^n)du.
\]
Then, by the definition of $\Omega_n$ and direct calculation,
\begin{align*}
&S_t^n=\lfloor nx\rfloor-\frac{\lambda}{n}\int_0^t \gamma(\mathcal{S}_u^n, \mathcal{I}_u^n)du+\mathcal{M}_t(f_1),\\
&I_t^n=\lfloor ny\rfloor+\frac{\lambda}{n}\int_0^t \gamma(\mathcal{S}_u^n, \mathcal{I}_u^n)du-\int_0^tI_u^ndu+\mathcal{M}_t(f_2), \notag
\end{align*}
where $\mathcal{S}_u^n=\{i:~\eta_u^n(i)=0\}$ and $\mathcal{I}_u^n=\{i:~\eta_u^n(i)=1\}$.
We define
\[
\widetilde{\mathcal{M}}_t(f_1)=\mathcal{M}_t(f_1)-<\mathcal{M}(f_1), \widehat{\mathcal{M}}(H_g)>_t
\]
and
\[
\widetilde{\mathcal{M}}_t(f_2)=\mathcal{M}_t(f_2)-<\mathcal{M}(f_2), \widehat{\mathcal{M}}(H_g)>_t,
\]
then, as we have recalled, $\{\widetilde{\mathcal{M}}_t(f_1)\}_{0\leq t\leq T_0}$ and $\{\widetilde{\mathcal{M}}_t(f_2)\}_{0\leq t\leq T_0}$ are both martingales under $\widehat{P}^{\omega,g}_{\lambda,n}$ and
\begin{align*}
&S_t^n=\lfloor nx\rfloor-\frac{\lambda}{n}\int_0^t \gamma(\mathcal{S}_u^n, \mathcal{I}_u^n)du+\widetilde{\mathcal{M}}_t(f_1)+<\mathcal{M}(f_1), \widehat{\mathcal{M}}(H_g)>_t,\\
&I_t^n=\lfloor ny\rfloor+\frac{\lambda}{n}\int_0^t \gamma(\mathcal{S}_u^n, \mathcal{I}_u^n)du-\int_0^tI_u^ndu+\widetilde{\mathcal{M}_t}(f_2)+<\mathcal{M}(f_2), \widehat{\mathcal{M}}(H_g)>_t. \notag
\end{align*}
By Equation \eqref{equ 3.1} and direct calculation,
\begin{align*}
&d<\mathcal{M}(f_1), \widehat{\mathcal{M}}(H_g)>_t=\frac{1}{H_g(t-,\eta_{t-}^n)}d<\mathcal{M}(f_1), \mathcal{M}(H_g)>_t\\
&=-\frac{H_g(t-,\eta_{t-}^n)}{H_g(t-,\eta_{t-}^n)}\big(e^{i_t(g)-s_t(g)}-1\big)\frac{\lambda}{n}\gamma(\mathcal{S}_t^n, \mathcal{I}_t^n)dt=-\big(e^{i_t(g)-s_t(g)}-1\big)\frac{\lambda}{n}\gamma(\mathcal{S}_t^n, \mathcal{I}_t^n)dt.
\end{align*}
According to a similar calculation,
\[
d<\mathcal{M}(f_1), \widehat{\mathcal{M}}(H_g)>_t=\big[-(e^{-i_t(g)}-1)I_t^n+(e^{i_t(g)-s_t(g)}-1)\frac{\lambda}{n}\gamma(\mathcal{S}_t^n, \mathcal{I}_t^n)\big]dt.
\]
As a result,
\begin{align}\label{equ StLDP}
\frac{S_t^n}{n}&=\frac{\lfloor nx\rfloor}{n}-\frac{\lambda}{n^2}\int_0^t e^{i_u(g)-s_u(g)}\gamma(\mathcal{S}_u^n, \mathcal{I}_u^n)du+\frac{1}{n}\widetilde{\mathcal{M}}_t(f_1)\\
&=\frac{\lfloor nx\rfloor}{n}-\int_0^t\lambda e^{i_u(g)-s_u(g)}\big(\frac{S_u^n}{n}\frac{I_u^n}{n}(E\rho)+\varepsilon_{u}^n\big) du+\frac{1}{n}\widetilde{\mathcal{M}}_t(f_1)\notag
\end{align}
and
\begin{align}\label{equ ItLDP}
\frac{I_t^n}{n}&=\frac{\lfloor ny\rfloor}{n}+\int_0^t-e^{-i_u(g)}\frac{I_u^n}{n}+\frac{\lambda}{n^2}e^{i_u(g)-s_u(g)}\gamma(\mathcal{S}_u^n, \mathcal{I}_u^n)du+\frac{1}{n}\widetilde{\mathcal{M}}_t(f_2)\\
&=\frac{\lfloor ny\rfloor}{n}+\int_0^t-e^{-i_u(g)}\frac{I_u^n}{n}+\lambda e^{i_u(g)-s_u(g)}\big(\frac{S_u^n}{n}\frac{I_u^n}{n}(E\rho)+\varepsilon_{u}^n\big)du+\frac{1}{n}\widetilde{\mathcal{M}}_t(f_2), \notag
\end{align}
where $\varepsilon_u^n=\frac{\gamma(\mathcal{S}_u^n, \mathcal{I}_u^n)-(E\rho)S_u^nI_u^n}{n^2}$.

As a result, to prove Lemma \ref{lemma 3.3 LLNunderChangeMeasure}, we only need to show that
$\sup_{0\leq u\leq T_0}|\varepsilon_u^n|$ converges in $\widehat{P}^g_{\lambda,n}$-probability to $0$ and $\sup_{0\leq t\leq T_0}|\frac{1}{n}\widetilde{\mathcal{M}}_t(f_i)|$ converges in $\widehat{P}^g_{\lambda,n}$-probability to $0$ for $i=1,2$.

Since $S_t^n, I_t^n\leq n$, there exists $K_2\in (0,+\infty)$ not depending on $n, \omega$ such $\Lambda_t^n(g)\leq e^{K_2n}$ for all $0\leq t\leq T_0$. Therefore,
\begin{equation}\label{equ 3.5}
\widehat{P}^g_{\lambda,n}\big(\sup_{0\leq u\leq T_0}|\varepsilon_u^n|\geq \epsilon\big)
\leq e^{K_2n}P_{\lambda,n}\big(\sup_{0\leq u\leq T_0}|\varepsilon_u^n|\geq \epsilon\big).
\end{equation}
Since $|\mathcal{S}_t^n|=S_t^n$ and $|\mathcal{I}_t^n|=I_t^n$, $\sup_{0\leq u\leq T_0}|\varepsilon_u^n|$ converges in $\widehat{P}^g_{\lambda,n}$-probability to $0$ according to Equation \eqref{equ 3.5} and Lemma \ref{lemma 3.1}.

To prove $\sup_{0\leq t\leq T_0}|\frac{1}{n}\widetilde{\mathcal{M}}_t(f_i)|$ converges in $\widehat{P}^g_{\lambda,n}$-probability to $0$ for $i=1,2$, we only need to show that
\[
[\frac{1}{n}\widetilde{\mathcal{M}}(f_i), \frac{1}{n}\widetilde{\mathcal{M}}(f_i)]_{T_0}=\frac{1}{n^2}[\widetilde{\mathcal{M}}(f_i), \widetilde{\mathcal{M}}(f_i)]_{T_0}
\]
converges in $\widehat{P}^g_{\lambda,n}$-probability to $0$ for $i=1,2$. As we have recalled, $[\widetilde{\mathcal{M}}(f_i), \widetilde{\mathcal{M}}(f_i)]=[\mathcal{M}(f_i), \mathcal{M}(f_i)]$. Since $S_t^n, I_t^n\leq n$, there exists $K_3, K_4\in (0,+\infty)$ not depending on $n,\omega$ such that $\{[\mathcal{M}(f_i), \mathcal{M}(f_i)]_t\}_{0\leq t\leq T_0}$ is stochastically dominated from above by $\{K_4Y_{nK_3t}\}_{0\leq t\leq T_0}$ under $P_{\lambda,n}^\omega$, where $\{Y_t\}_{t\geq 0}$ is a Poison process with rate $1$. Then, for any $\epsilon>0$,
\[
P_{\lambda,n}\big(\frac{1}{n^2}[\mathcal{M}(f_i), \mathcal{M}(f_i)]_{T_0}\geq \epsilon\big)\leq e^{-n^2\epsilon}Ee^{K_4Y_{nK_3T_0}}=e^{-n^2\epsilon}e^{(e^{K_4}-1)nK_3T_0}
\]
and hence
\[
\widehat{P}^g_{\lambda,n}\big(\frac{1}{n^2}[\widetilde{\mathcal{M}}(f_i), \widetilde{\mathcal{M}}(f_i)]_{T_0}\geq \epsilon\big)\leq e^{-n^2\epsilon}e^{(e^{K_4}-1)nK_3T_0}e^{K_2n}.
\]
Therefore, $[\frac{1}{n}\widetilde{\mathcal{M}}(f_i), \frac{1}{n}\widetilde{\mathcal{M}}(f_i)]_{T_0}$ converges in $\widehat{P}^g_{\lambda,n}$-probability to $0$ for $i=1,2$ and the proof is complete.

\qed

At the end of this section, we give the proof of Equation \eqref{equ LDP lower}.

\proof[Proof of Equation \eqref{equ LDP lower}]

If $\inf_{O\bigcap \mathcal{B}}(I_{ini}(f_0)+I_{dyn}(f))=+\infty$, then the conclusion is trivial, so we only deal with the case where
$\inf_{O\bigcap \mathcal{B}}(I_{ini}(f_0)+I_{dyn}(f))<+\infty$. For any $\epsilon>0$, there exists $f^\epsilon\in O\bigcap \mathcal{B}$ such that
\[
I_{ini}(f^\epsilon_0)+I_{dyn}(f^\epsilon)<\inf_{O\bigcap \mathcal{B}}(I_{ini}(f_0)+I_{dyn}(f))+\epsilon.
\]
Hence, by Theorem \ref{theorem 2.2 Idyn Expression}, $s_0(f^\epsilon)\geq 0, i_0(f^\epsilon)\geq 0$ and $s_0(f^\epsilon)+i_0(f^\epsilon)\leq 1$ while $f^\epsilon$ is absolutely continuous. As we have shown in the proof of Theorem \ref{theorem 2.2 Idyn Expression}, there exists $h_1^\epsilon(t), h_2^\epsilon(t)\geq 0$ such that
\[
\begin{cases}
s_t^\prime(f^\epsilon)&=-h_2^\epsilon(t)\lambda(E\rho)i_t(f^\epsilon)s_t(f^\epsilon),\\
i_t^\prime(f^\epsilon)&=-h_1^\epsilon(t)i_t(f^\epsilon)+h_2^\epsilon(t)\lambda(E\rho)i_t(f^\epsilon)s_t(f^\epsilon).
\end{cases}
\]
Similarly with that in the proof of Theorem \ref{theorem 2.2 Idyn Expression}, we let $h_1^n\in C^2\big([0,T_0], \mathbb{R}\big)$ such that $h_1^n(t)>0$ for all $t\in [0, T_0]$ while
\[
\lim_{n\rightarrow+\infty}\int_0^{T_0}|h_1^n(t)-h_1^\epsilon(t)|dt=0
\]
and
$h_2^n\in C^2\big([0,T_0], \mathbb{R}\big)$ such that $h_2^n(t)>0$ for all $t\in [0, T_0]$ while
\[
\lim_{n\rightarrow+\infty}\int_0^{T_0}|h_2^n(t)-h_2^\epsilon(t)|dt=0
\]
and define $g^n\in C^2\big([0,T_0], \mathbb{R}^2\big)$ such that $s_t(g^n)$ and $i_t(g^n)$ satisfies
\[
h_1^n(t)=\exp\big\{-i_t(g^n)\big\} \text{~and~} h_2^n(t)=\exp\big\{i_t(g^n)-s_t(g^n)\big\}
\]
for all $t\in [0, T_0]$. Then we define $f^n\in C^2\big([0,T_0], \mathbb{R}^2\big)$ as the solution to the ODE
\begin{align*}
\begin{cases}
&s_t^\prime(f^n)=-h_2^n(t)\lambda(E\rho)i_t(f^n)s_t(f^n),\\
&i_t^\prime(f^n)=-h_1^n(t)i_t(f^n)+h_2^n(t)\lambda(E\rho)i_t(f^n)s_t(f^n),\\
&(s_0(f^n), i_0(f^n))=(s_0(f^\epsilon), i_0(f^\epsilon)).
\end{cases}
\end{align*}
Then, according to our assumptions of $h_1^n, h_2^n$ and Grownwall's inequality, $(f^n)^\prime$ converges to $(f^\epsilon)^\prime$ in $L_1\big([0,T_0], \mathbb{R}^2\big)$ while $f^n$ converges to $f^\epsilon$ uniformly on $[0, T_0]$. Note that the definition of $f^n$ ensures that $f^n\in \mathcal{B}$. Since $g^n\in  C^2\big([0,T_0], \mathbb{R}^2\big)$, according to a similar analysis with that in the proof of Theorem \ref{theorem 2.2 Idyn Expression},
\begin{align*}
I_{dyn}(f^n)&=\Phi_{f^n}(g^n)\\
&=\int_0^{T_0}L_t(f^n)+(i_t(f^n)+s_t(f^n))^\prime+i_t(f^n)+s_t^\prime(f^n)+\lambda(E\rho)i_t(f^n)s_t(f^n) dt
\end{align*}
and hence
\[
\lim_{n\rightarrow +\infty} I_{dyn}(f^n)=I_{dyn}(f^\epsilon).
\]
Therefore, there exists $n_1\geq 1$ such that $|I_{dyn}(f^{n_1})-I_{dyn}(f)|<\epsilon$ and $f^{n_1}\in O\bigcap B$.

According to the definition of $g^{n_1}$ and $f^{n_1}$,
\begin{equation*}
\begin{cases}
&\frac{d}{dt}s_t(f^{n_1})=-e^{i_t(g^{n_1})-s_t(g^{n_1})}\lambda(E\rho)s_t(f^{n_1})i_t(f^{n_1}),\\
&\frac{d}{dt}i_t(f^{n_1})=-e^{-i_t(g^{n_1})}i_t(g^{n_1})+e^{i_t(g^{n_1})-s_t(g^{n_1})}\lambda(E\rho)s_t(f^{n_1})i_t(f^{n_1}), \\
&(s_0(f^{n_1}),i_0(f^{n_1}))=(s_0(f^\epsilon),i_0(f^\epsilon)).
\end{cases}
\end{equation*}
Therefore, according to Lemma \ref{lemma 3.3 LLNunderChangeMeasure}, $\vartheta^n$ converges in $\widehat{P}^{g^{n_1},s_0(f^\epsilon),i_0(f^\epsilon)}_{\lambda,n}$-probability to $f^{n_1}$ as $n\rightarrow+\infty$.

By direct calculation,
\begin{align*}
\Lambda^n_{T_0}(g^{n_1})=&\exp\Big\{n\vartheta^n_{T_0}\cdot g^{n_1}_{T_0}-n\vartheta^n_0\cdot g^{n_1}_0-\int_0^{T_0}n\vartheta^n_t\cdot (g_t^{n_1})^\prime dt\\
&-\int_0^{T_0}\big(e^{g_t^{n_1}\cdot l_1}-1\big)I_t^n+\big(e^{g_t^{n_1}\cdot l_2}-1\big)\frac{\lambda}{n}\gamma(\mathcal{S}_t^n, \mathcal{I}_t^n)dt\Big\}.
\end{align*}
Let $\delta_n$ be defined as before Lemma \ref{lemma 3.1}. According to the above expression of $\Lambda^n_{T_0}(g^{n_1})$, for given $\epsilon>0$, there exists $\delta_0>0$ not depending on $n$ such that
\[
\Lambda_{T_0}^n(g^{n_1})\leq \exp\big\{n(\Phi_{f^{n_1}}(g^{n_1})+\epsilon+\frac{\|g^{n_1}\|}{n})\big\}
\]
conditioned on $\vartheta^n\in B(f^{n_1}, \delta_0)$ and $\delta_n\leq \delta_0$, where $B(f^{n_1}, r)$ is the ball concentrated on $f^{n_1}$ with radius $r$. Since $O$ is open, we can further assume that $\delta_0$ makes $B(f^{n_1}, \delta_0)\subseteq O$. As a result,
\begin{align*}
&P_{\lambda,n}\big(\vartheta^n\in O\big)\geq P_{\lambda,n}\big(\vartheta^n\in B(f^{n_1},\delta_0), \delta_n\leq \delta_0, S_0^n=\lfloor ns_0(f^\epsilon)\rfloor, I_0^n=\lfloor ni_0(f^\epsilon)\rfloor\big) \\
&=\widehat{E}_{\lambda,n}^{g^{n_1}}\Big[\frac{1}{\Lambda_{T_0}^n(g^{n_1})}1_{\{\vartheta^n\in B(f^{n_1},\delta_0), \delta_n\leq \delta_0, S_0^n=\lfloor ns_0(f^\epsilon)\rfloor, I_0^n=\lfloor ni_0(f^\epsilon)\rfloor\}}\Big] \\
&\geq \exp\big\{-n(\Phi_{f^{n_1}}(g^{n_1})+\epsilon+\frac{\|g^{n_1}\|}{n})\big\}\widehat{P}_{\lambda,n}^{g^{n_1}, s_0(f^\epsilon), i_0(f^\epsilon)}\big(\vartheta^n\in B(f^{n_1},\delta_0), \delta_n\leq \delta_0\big) \\
&\times \widehat{P}_{\lambda,n}^{g^{n_1}}\big(S_0^n=\lfloor ns_0(f^\epsilon)\rfloor, I_0^n=\lfloor ni_0(f^\epsilon)\rfloor\big).
\end{align*}
As we have shown, $\vartheta^n$ converges in $\widehat{P}^{g^{n_1},s_0(f^\epsilon),i_0(f^\epsilon)}_{\lambda,n}$-probability to $f^{n_1}$ as $n\rightarrow+\infty$. Further, according to the analysis in the proof of Lemma \ref{lemma 3.3 LLNunderChangeMeasure}, $\delta_n$ converges in $\widehat{P}^{g^{n_1},s_0(f^\epsilon),i_0(f^\epsilon)}_{\lambda,n}$-probability to $0$ as $n\rightarrow+\infty$. Therefore,
\[
\lim_{n\rightarrow+\infty}\widehat{P}_{\lambda,n}^{g^{n_1}, s_0(f^\epsilon), i_0(f^\epsilon)}\big(\vartheta^n\in B(f^{n_1},\delta_0), \delta_n\leq \delta_0\big)=1
\]
and
\begin{align*}
&\liminf_{n\rightarrow+\infty}\frac{1}{n}\log P_{\lambda,n}\big(\vartheta^n\in O\big)\geq \\
&-\Phi_{f^{n_1}}(g^{n_1})-\epsilon+\lim_{n\rightarrow+\infty}\frac{1}{n}\log \widehat{P}_{\lambda,n}^{g^{n_1}}\big(S_0^n=\lfloor ns_0(f^\epsilon)\rfloor, I_0^n=\lfloor ni_0(f^\epsilon)\rfloor\big).
\end{align*}
Since $\{\Lambda_t^{n}(g^{n_1})\}_{0\leq t\leq T_0}$ is an exponential martingale with $\Lambda_0^{n}(g^{n_1})=1$, $\vartheta^n_0$ have the same distribution under $\widehat{P}_{\lambda,n}^{g^{n_1}}$ and $P_{\lambda,n}$. As a result, according to Assumption \textbf{A}, Theorem \ref{theorem 2.2 Idyn Expression} and Strling's formula,
\begin{align*}
&\lim_{n\rightarrow+\infty}\frac{1}{n}\log \widehat{P}_{\lambda,n}^{g^{n_1}}\big(S_0^n=\lfloor ns_0(f^\epsilon)\rfloor, I_0^n=\lfloor ni_0(f^\epsilon)\rfloor\big)\\
&=\lim_{n\rightarrow+\infty}\frac{1}{n}\log P_{\lambda,n}\big(S_0^n=\lfloor ns_0(f^\epsilon)\rfloor, I_0^n=\lfloor ni_0(f^\epsilon)\rfloor\big) \\
&=\lim_{n\rightarrow+\infty}\frac{1}{n}\log \Bigg\{{n \choose \lfloor ns_0(f^\epsilon)\rfloor}{n-\lfloor ns_0(f^\epsilon)\rfloor \choose \lfloor ni_0(f^\epsilon)\rfloor}p_0^{\lfloor ns_0(f^\epsilon)\rfloor}\\
&\text{~\quad\quad}\times p_1^{\lfloor ni_0(f^\epsilon)\rfloor}(1-p_0-p_1)^{n-\lfloor ns_0(f^\epsilon)\rfloor-\lfloor ni_0(f^\epsilon)\rfloor}\Bigg\}\\
&=-\Bigg(s_0(f^\epsilon)\log \frac{s_0(f^\epsilon)}{p_0}+i_0(f^\epsilon)\log \frac{i_0(f^\epsilon)}{p_1}\\
&\text{\quad\quad}+(1-s_0(f^\epsilon)-i_0(f^\epsilon))\log\frac{1-s_0(f^\epsilon)-i_0(f^\epsilon)}{1-p_0-p_1}\Bigg)=-I_{ini}(f^\epsilon_0).
\end{align*}
Therefore,
\begin{align*}
&\liminf_{n\rightarrow+\infty}\frac{1}{n}\log P_{\lambda,n}\big(\vartheta^n\in O\big)\geq -\Phi_{f^{n_1}}(g^{n_1})-\epsilon-I_{ini}(f^\epsilon_0) \\
&=-I_{dyn}(f^{n_1})-I_{ini}(f^\epsilon_0)-\epsilon\geq -I_{dyn}(f^\epsilon)-I_{ini}(f^\epsilon_0)-2\epsilon\\
&\geq -\inf_{O\bigcap \mathcal{B}}(I_{ini}(f_0)+I_{dyn}(f))-3\epsilon.
\end{align*}
Since $\epsilon$ is arbitrary, the proof is complete.

\qed

\proof[Remark]

We think that the strategy of the above proof of Equation \eqref{equ LDP lower} can be utilized in the study of large deviations for some other density-dependent Markov processes. Roughly speaking, the core idea of the proof is to show that $\vartheta^n$ converges to $f$ under $\widehat{P}^g$ and consequently the rate function $I(f)$ satisfies $nI(f)\approx -\log \frac{dP}{d\widehat{P}^g}\Big|_{\vartheta^n=f}$, where
\[
f^\prime=l_1e^{g\cdot l_1}H_1(f)+l_2e^{g\cdot l_2}H_2(f)
\]
while $nH_i(f_t)$ is nearly the rate at which $\vartheta^n$ flips from $nf_t$ to $nf_t+l_i$.
Similarly, for a density-dependent Markov process $\{X_t^n\}_{t\geq 0}$ with parameters
$\{F_l\}_{l\in \mathcal{A}}$, let
\begin{align*}
I(f)=&\sup_{g\in C^2\big([0,T_0],\mathbb{R}^d\big)}\Big\{f_{T_0}\cdot g_{T_0}-f_0\cdot g_0-\int_0^{T_0}f_t\cdot g_t^\prime dt\\
&-\int_0^{T_0}\sum_{l\in \mathcal{A}}\big(e^{g_t\cdot l}-1\big)F_l(f_t)dt\Big\}.
\end{align*}
If one could show that $I(f)<+\infty$ implies that there exists $g$ such that
\[
f^\prime=\sum_{l\in \mathcal{A}}le^{g\cdot l}F_l(f),
\]
which is an analogue of Theorem \ref{theorem 2.2 Idyn Expression} and intuitively holds according to a non-rigorous variational method, then the large deviation principle of $\{\frac{X_t^n}{n}\}_{0\leq t\leq T_0}$ with rate function $I$ would hold according to the above strategy. We guess that this analysis may work for all the cases where $\mathcal{A}$ is finite and $\{F_l\}_{l\in \mathcal{A}}$ are bounded and smooth. However, we have not yet found a rigorous proof of the above analogue of Theorem \ref{theorem 2.2 Idyn Expression} for these general cases. We will work on this question as a further investigation.

\qed

\section{The proof of Equation \eqref{equ LDP upper}}\label{section four}
In this section we give the proof of Equation \eqref{equ LDP upper}. First we show that this equation holds for compact sets.
\begin{lemma}\label{lemma 4.1}
For any compact set $C\subseteq \mathcal{D}\big([0,T_0], \mathbb{R}^2\big)$,
\begin{equation*}
\limsup_{n\rightarrow+\infty}\frac{1}{n}\log P_{\lambda, n}\big(\vartheta^n\in C\big)\leq -\inf_{f\in C\bigcap \mathcal{B}}\big(I_{dyn}(f)+I_{ini}(f_0)\big).
\end{equation*}
\end{lemma}

\proof

Let $\beta_{2,x}(y)=y\cdot x-\log\big(1-p_0-p_1+e^{s_y}p_0+e^{i_y}p_1$, then $I_{ini}(x)=\sup_{y\in \mathbb{R}^2}\beta_{2,x}(y)$. For any $\epsilon>0$, by Lemma \ref{lemma 3.1} and the fact that $\vartheta\in \mathcal{B}$ almost surely,
\begin{equation}\label{equ 4.1}
\limsup_{n\rightarrow+\infty}\frac{1}{n}\log P_{\lambda, n}\big(\vartheta^n\in C\big)=\limsup_{n\rightarrow+\infty}\frac{1}{n}\log P_{\lambda, n}\big(\vartheta^n\in C\bigcap \mathcal{B},\delta_n\leq \epsilon\big).
\end{equation}
For any $g\in C^2\big([0,T_0], \mathbb{R}^2\big)$ and $y\in \mathbb{R}^2$, by the expression of $\Lambda_{T_0}^n(g)$ given in Section \ref{section three}, conditioned on $\vartheta^n\in C\bigcap \mathcal{B}$ and $\delta_n\leq \epsilon$,
\[
e^{ny\cdot \vartheta_0^n}\Lambda^n_{T_0}(g)\geq \exp\big\{n\big[\inf_{f\in C\bigcap \mathcal{B}}\big(y\cdot f_0+\Phi_f(g)\big)-\lambda\epsilon e^{\|g\|+1}\big]\big\}.
\]
Hence, by Lemma \ref{lemma 3.2 LambdaIsMartingale}, for any $g\in C^2\big([0,T_0], \mathbb{R}^2\big)$ and $y\in \mathbb{R}^2$,
\begin{align*}
E_{\lambda,n}e^{ny\cdot\vartheta_0^n}&=E_{\lambda,n}\big(e^{ny\cdot\vartheta_0^n}\Lambda_{0}^{n}(g)\big)\\
&=E_{\lambda,n}\big(e^{ny\cdot\vartheta_0^n}\Lambda_{T_0}^{n}\big)\geq E_{\lambda,n}\big(e^{ny\cdot\vartheta_0^n}\Lambda_{T_0}^{n}1_{\{\vartheta^n\in C\bigcap \mathcal{B},\delta_n\leq \epsilon\}}\big)\\
&\geq \exp\Big\{n\big(\inf_{f\in C\bigcap \mathcal{B}}\big(y\cdot f_0+\Phi_f(g)\big)-\lambda\epsilon e^{\|g\|+1}\big)\Big\}\\
&\text{~\quad\quad}\times P_{\lambda,n}\big(\vartheta^n\in C\bigcap \mathcal{B},\delta_n\leq \epsilon\big).
\end{align*}
By Assumption \textbf{A}, $E_{\lambda,n}e^{ny\cdot\vartheta_0^n}=e^{n\log\big(1-p_0-p_1+e^{s_y}p_0+e^{i_y}p_1\big)}$, hence,
\begin{align*}
&\limsup_{n\rightarrow+\infty}\frac{1}{n}\log P_{\lambda,n}\big(\vartheta^n\in C\bigcap \mathcal{B},\delta_n\leq \epsilon\big)\\
&\leq -\inf_{f\in C\bigcap \mathcal{B}}\big(y\cdot f_0+\Phi_f(g)\big)+\log\big(1-p_0-p_1+e^{s_y}p_0+e^{i_y}p_1\big)+\lambda\epsilon e^{\|g\|+1}\\
&=-\inf_{f\in C\bigcap \mathcal{B}}\big(\beta_{2,f_0}(y)+\Phi_f(g)\big)+\lambda\epsilon e^{\|g\|+1}.
\end{align*}
By Equation \eqref{equ 4.1},
\[
\limsup_{n\rightarrow+\infty}\frac{1}{n}\log P_{\lambda, n}\big(\vartheta^n\in C\big)\leq -\inf_{f\in C\bigcap \mathcal{B}}\big(\beta_{2,f_0}(y)+\Phi_f(g)\big)+\lambda\epsilon e^{\|g\|+1}
\]
and hence
\[
\limsup_{n\rightarrow+\infty}\frac{1}{n}\log P_{\lambda, n}\big(\vartheta^n\in C\big)\leq -\inf_{f\in C\bigcap \mathcal{B}}\big(\beta_{2,f_0}(y)+\Phi_f(g)\big)
\]
since $\epsilon$ is arbitrary. Since $g$ and $y$ are arbitrary,
\begin{equation}\label{equ 4.2}
\limsup_{n\rightarrow+\infty}\frac{1}{n}\log P_{\lambda, n}\big(\vartheta^n\in C\big)\leq -\sup_{y,g}\inf_{f\in C\bigcap \mathcal{B}}\big(\beta_{2,f_0}(y)+\Phi_f(g)\big).
\end{equation}
Since $\beta_{2,f_0}(y)+\Phi_f(g)$ is concave with $(g,y)$ while convex with $f$ and $C$ is compact, according to the minimax theorem given in \cite{Sion1958},
\begin{equation}\label{equ 4.3}
\sup_{y,g}\inf_{f\in C\bigcap \mathcal{B}}\big(\beta_{2,f_0}(y)+\Phi_f(g)\big)
=\inf_{f\in C\bigcap \mathcal{B}}\sup_{y,g}\big(\beta_{2,f_0}(y)+\Phi_f(g)\big).
\end{equation}
For given $f$,
\[
\sup_{y,g}\big(\beta_{2,f_0}(y)+\Phi_f(g)\big)=\sup_{y}\beta_{2,f_0}(y)+\sup_{g}\Phi_f(g)=I_{ini}(f_0)+I_{dyn}(f)
\]
and hence Lemma \ref{lemma 4.1} follows from Equations \eqref{equ 4.2} and \eqref{equ 4.3}.

\qed

At the end of this section, we give the proof of Equation \eqref{equ LDP upper}.

\proof[Proof of Equation \eqref{equ LDP upper}]

By Lemma \ref{lemma 4.1}, we only need to show that $\{\vartheta^n\}_{n\geq 1}$ is exponential tight, which is equivalent to the following two properties
(see the main theorem of \cite{Puhalskii1994}).

(\textbf{1}) \[
\limsup_{M\rightarrow+\infty}\limsup_{n\rightarrow+\infty}\frac{1}{n}\log P_{\lambda,n}\big(\|\vartheta^n\|>M\big)=-\infty.
\]
(\textbf{2}) For any $\epsilon>0$,
\[
\limsup_{\delta\rightarrow 0}\limsup_{n\rightarrow+\infty}\frac{1}{n}\log \sup_{\tau \in \mathcal{T}_0}
P_{\lambda,n}\big(\sup_{0\leq t\leq \delta}|\frac{I_{t+\tau}^n-I_\tau^n}{n}|+|\frac{S_{t+\tau}^n-S_\tau^n}{n}|>\epsilon\big)=-\infty,
\]
where $\mathcal{T}_0$ is the set of stopping times of $\{\eta_t^n\}_{0\leq t\leq T_0}$ with upper bound $T_0$.

To check Property (\textbf{1}), note that $S_t^n+I_t^n\leq n$ implies that
\[
P_{\lambda,n}\big(\|\vartheta^n\|>1\big)=0,
\]
Property (\textbf{1}) follows from which directly.

Now we only need to check Property (\textbf{2}). By Lemma \ref{lemma 3.1}, we only need to check
(\textbf{\.{2}}):
\[
\limsup_{\delta\rightarrow 0}\limsup_{n\rightarrow+\infty}\frac{1}{n}\log \sup_{\tau \in \mathcal{T}_0}
P_{\lambda,n}\big(\sup_{0\leq t\leq \delta}|\frac{I_{t+\tau}^n-I_\tau^n}{n}|+|\frac{S_{t+\tau}^n-S_\tau^n}{n}|>\epsilon, \delta_n\leq M\big)=-\infty
\]
for any $M,\epsilon>0$. Conditioned on $\delta_n\leq M$, $S_t^n ,I_t^n\leq n$ implies that $\{|I_{t+\tau}^n-I_\tau^n|+|S_{t+\tau}^n-S_\tau^n|\}_{t\geq 0}$ is stochastically dominated from above by $\{2Y_{K_5nt}\}_{t\geq 0}$ for some $K_5=K_5(M)\in (0,+\infty)$ not depending on $n$, where $\{Y_t\}_{t\geq 0}$ is the Poisson process with rate $1$. Therefore, by Chebyshev's inequality, for any $\theta>0$,
\begin{align*}
&P_{\lambda,n}\big(\sup_{0\leq t\leq \delta}|\frac{I_{t+\tau}^n-I_\tau^n}{n}|+|\frac{S_{t+\tau}^n-S_\tau^n}{n}|>\epsilon, \delta_n\leq M\big) \\
&\leq P(e^{\theta Y_{K_5n\delta}}>e^{\frac{\theta n\epsilon}{2}})=e^{-\frac{\theta n\epsilon}{2}}e^{K_5n\delta(e^\theta-1)}
\end{align*}
and then,
\[
\limsup_{\delta\rightarrow 0}\limsup_{n\rightarrow+\infty}\frac{1}{n}\log \sup_{\tau \in \mathcal{T}_0}
P_{\lambda,n}\big(\sup_{0\leq t\leq \delta}|\frac{I_{t+\tau}^n-I_\tau^n}{n}|+|\frac{S_{t+\tau}^n-S_\tau^n}{n}|>\epsilon, \delta_n\leq M\big)\leq -\frac{\theta\epsilon}{2},
\]
Property (\textbf{\.{2}}) follows from which since $\theta$ is arbitrary and hence Property (\textbf{2}) holds.

Since Properties (\textbf{1}) and (\textbf{2}) hold, $\{\vartheta^n\}_{n\geq 1}$ is exponential tight. Consequently, Equation \eqref{equ LDP upper} follows from Lemma \ref{lemma 4.1}.

\qed

\section{The outline of proof of Theorem \ref{Theorem 2.4 MDPofSIR}}\label{section five}
In this section, we give the outline of the proof of Theorem \ref{Theorem 2.4 MDPofSIR}. For later use, we first prove Lemma \ref{lemma 2.3 LLNofSIR} and Theorem \ref{theorem 2.5 JdynEpress}.

\proof[Proof of Lemma \ref{lemma 2.3 LLNofSIR}]

For $f\in \mathcal{B}$, we claim that $I_{ini}(f_0)+I_{dyn}(f)=0$ if and only if $f=\widehat{x}$, which we will prove at the end of this proof. Consequently,
$\inf\{I_{ini}(f_0)+I_{dyn}(f):~f\in\mathcal{B}\text{~and~}\|f-\widehat{x}\|\geq \epsilon\}>0$ for any $\epsilon>0$ and Lemma \ref{lemma 2.3 LLNofSIR} follows from Equation \eqref{equ LDP upper} directly.

Hence, we only need to prove our claim. By Theorem \ref{theorem 2.2 Idyn Expression}, $I_{ini}(x)$ gets minimum $0$ when and only when $x$ is the solution to
$\frac{\partial}{\partial s_x}I_{ini}(x)=\frac{\partial}{\partial i_x}I_{ini}(x)=0$, i.e, $s_x=p_0, i_x=p_1$. Then, for $f\in \mathcal{B}$ satisfying $I_{ini}(f_0)+I_{dyn}(f)=0$, i.e., $I_{ini}(f_0)=I_{dyn}(f)=0$, we have $s_0(f)=p_0$ and $i_0(f)=p_1$. For such $f$, since $I_{dyn}(f)=0$,
\[
\sup\{\Phi_f(cg):~c\in \mathbb{R}\}=\Phi_f(0g)=0
\]
for any given $g\in C^2\big([0,T_0], \mathbb{R}^2\big)$. Therefore, $\frac{d}{dc}\Phi_f(cg)\Big|_{c=0}=0$, i.e,
\[
\int_0^{T_0} \big(f_t^\prime-l_1i_t(f)-l_2\lambda(E\rho)s_t(f)i_t(f)\big)\cdot g_tdt=0
\]
for any $g\in C^2\big([0,T_0], \mathbb{R}^2\big)$, where $l_1=(0,-1)^T, l_2=(-1,1)^T$ defined as in Section \ref{section two}. Therefore, $f_t^\prime=l_1i_t(f)+l_2\lambda(E\rho)s_t(f)i_t(f)$ for $0\leq t\leq T_0$, i.e.,
\[
\begin{cases}
&s_t^\prime(f)=-\lambda(E\rho)s_t(f)i_t(f),\\
&i_t^\prime(f)=-i_t(f)+\lambda(E\rho)s_t(f)i_t(f),
\end{cases}
\]
$f=\widehat{x}$ follows from which and the fact that $s_0(f)=p_0, i_0(f)=p_1$.

\qed

\proof[Proof of Theorem \ref{theorem 2.5 JdynEpress}]

For given $x\in \mathbb{R}^2, y\in \mathbb{R}^2\setminus \{(0,0)^T\}$ and any $c\in R$,
\[
(cy)\cdot x-\frac{1}{2}(cy)^TM_{_0}(cy)=c(y\cdot x)-\frac{c^2}{2}y^TM_{_0}y
\]
and hence gets maximum $\frac{(y\cdot x)^2}{2y^TM_{_0}y}$ when $c=\frac{y\cdot x}{y^TM_{_0}y}$. Therefore,
\[
J_{ini}(x)=\sup\{\frac{(y\cdot x)^2}{2y^TM_{_0}y}:~y\neq (0,0)^T\}.
\]
For any $y\in \mathbb{R}^2$, by Cauchy-Schwartz's inequality,
\[
|y\cdot x|^2=|(M_{_0}y)\cdot (M_{_0}^{-\frac{1}{2}}x)|^2\leq y^TM_{_0}yx^TM^{-1}_{_0}x
\]
and hence $J_{ini}(x)\leq \frac{1}{2}x^TM^{-1}_{_0}x$. On the other hand, let $y_{_0}=M_0^{-1}x$, then
\[
J_{ini}(x)\geq \frac{(y_0\cdot x)^2}{2y_0^TM_{_0}y_0}=\frac{1}{2}x^TM^{-1}_{_0}x.
\]
Therefore, $J_{ini}(x)=\frac{1}{2}x^TM^{-1}_{_0}x$.

For $f\in \mathcal{D}\big([0,T_0], \mathbb{R}^2\big)$ satisfying $J_{dyn}(f)<+\infty$,
\[
J_{dyn}(f)=\frac{1}{2}\int_0^{T_0}(f_t^\prime-b_tf_t)^T\sigma_t^{-1}(f_t^\prime-b_tf_t) dt
\]
 holds as a special example of Equation (2.2) of Reference \cite{Xue2019}, the proof of which utilizes Riesz's representation Theorem to show that $J_{dyn}(f)<+\infty$ implies that $f$ is absolutely continuous and there exists $\varphi$ such that
\[
f^\prime_t-b_tf_t=\sigma_t\varphi_t
\]
while $J_{dyn}(f)=\frac{1}{2}\int_0^{T_0}\varphi_t^T\sigma_t\varphi_tdt$. Details of this proof could be checked in \cite{Xue2019}, which we omit here.

\qed

As another preparation work, we need the following lemma, which is an analogue of Lemma \ref{lemma 3.1}.

\begin{lemma}\label{lemma 5.1}
For any $\epsilon>0$,
\[
\lim_{n\rightarrow+\infty}\frac{n}{a_n^2}\log P(\frac{n\delta_n}{a_n}>\epsilon)=-\infty.
\]
\end{lemma}

\proof[Proof of Lemma \ref{lemma 5.1}]

According to Chebyshev's's inequality and similar analysis with that in the proof of Lemma \ref{lemma 3.1}, for any $C,D\in C_n$ such that $C\bigcap D=\emptyset$ and $\theta>0$,
\[
P\big(\gamma(C,D)-|C||D|(E\rho)>na_n\epsilon\big)\leq e^{-\theta a_n^2\epsilon}\big[Ee^{\frac{a_n}{n}\theta(\rho-E\rho)}\big]^{n^2}.
\]
Since $\frac{a_n}{n}\rightarrow 0$, by Taylor's expansion formula,
\begin{align*}
Ee^{\frac{a_n}{n}\theta(\rho-E\rho)}&=1+\frac{a_n}{n}\theta E(\rho-E\rho)+\frac{a_n^2}{2n^2}\theta^2{\rm Var}(\rho)+o(\frac{a_n^2}{n^2})\\
&=1+\frac{a_n^2}{2n^2}\theta^2{\rm Var}(\rho)+o(\frac{a_n^2}{n^2})
\end{align*}
and hence
\[
P\big(\gamma(C,D)-|C||D|(E\rho)>na_n\epsilon\big)\leq e^{-\theta a_n^2\epsilon+\frac{a_n^2}{2}[{\rm Var}(\rho)\theta^2+o(1)]}
\]
according to the fact that $1+x\leq e^x$. Let $\theta=\frac{\epsilon}{{\rm Var}(\rho)}$, then
\[
P\big(\gamma(C,D)-|C||D|(E\rho)>na_n\epsilon\big)\leq e^{-\frac{a_n^2\epsilon^2}{2{\rm Var}(\rho)}[1+o(1)]}.
\]
Note that $o(1)$ in the above inequality does not rely on $C,D$. According to a similar analysis,
\[
P\big(\gamma(C,D)-|C||D|(E\rho)<-na_n\epsilon\big)\leq e^{-\frac{a_n^2\epsilon^2}{2{\rm Var}(\rho)}[1+o(1)]}.
\]
Then, since the number of subsets of $C_n$ is $2^n$,
\[
P(\frac{n\delta_n}{a_n}>\epsilon)\leq e^{4n}e^{-\frac{a_n^2\epsilon^2}{2{\rm Var}(\rho)}[1+o(1)]},
\]
Lemma \ref{lemma 5.1} follows from which directly since $\frac{n}{a_n^2}\rightarrow 0$.

\qed

Similarly with that in Section \ref{section three}, for any $g\in C^2\big([0,T_0], \mathbb{R}^2\big)$, we define
\[
\zeta_g(t,\eta)=\frac{a_n}{n}\Big[i_t(g)\big(\sum_{i=1}^n1_{\{\eta(i)=1\}}-n\widehat{i}_t\big)+s_t(g)\big(\sum_{i=1}^n1_{\{\eta(i)=0\}}-n\widehat{s}_t\big)\Big],
\]
where $\widehat{x}_t=(\widehat{s}_t,\widehat{i}_t)^T$ defined as in Section \ref{section two}. Then,
\[
\zeta_g(t,\eta_t^n)=\frac{a_n}{n}g_t\cdot\big[(S_t^n,I_t^n)^T-n\widehat{x}_t\big].
\]
We further define $V_g(t,\eta_t^n)=e^{\zeta_g(t,\eta_t^n)}$ and
\[
\Xi_t^n(g)=\frac{V_g(t,\eta_t^n)}{V_g(0,\eta_0^n)}\exp\Big(-\int_0^t\frac{(\frac{\partial}{\partial u}+\Omega_n)V_g(u,\eta_u^n)}{V_g(u,\eta^n_u)}du\Big),
\]
then we have the following lemma, which is an analogue of Lemma \ref{lemma 3.2 LambdaIsMartingale}.

\begin{lemma}\label{lemma 5.2 XiIsMartingale}
For any $\omega\in X$ and $g\in C^2\big([0,T_0], \mathbb{R}^2\big)$, $\{\Xi_t^n(g)\}_{0\leq t\leq T_0}$ is a martingale with expectation $1$ under the quenched measure
$P_{\lambda,n}^\omega$.
\end{lemma}
The proof of Lemma \ref{lemma 5.2 XiIsMartingale} is nearly the same as that of Lemma \ref{lemma 3.2 LambdaIsMartingale}, which we omit.

According to the definition of $\Omega_n$ and Taylor's expansion formula,
\begin{align}\label{equ 5.1 ExpressionXi}
\Xi_t^n(g)=\exp\Bigg\{&\frac{a_n^2}{n}\Bigg(g_t\cdot \nu_t^n-g_0\cdot \nu_0^n-\int_0^tg_u^\prime\cdot \nu_u^n+(b_u^n\nu_n^n)\cdot g_u^n \notag\\
&+\lambda g_u\cdot l_2(\frac{n\varepsilon_u^n}{a_n})+\frac{1}{2}g_u^T\sigma_u^n g_u+\frac{\lambda}{2}g_u^T(l_2\varepsilon_u^nl_2^T)g_u du+o(1)\Bigg)\Bigg\},
\end{align}
where $\varepsilon_u^n=\frac{\gamma(\mathcal{S}_u^n, \mathcal{I}_u^n)-(E\rho)S_u^nI_u^n}{n^2}$ defined as in Section \ref{section three},
\[
\sigma_u^n=\sum_{i=1}^2l_iH_i(\vartheta_t^n)l_i^T=
\begin{pmatrix}
\lambda(E\rho)\frac{I_u^n}{n}\frac{S_u^n}{n} & -\lambda(E\rho)\frac{I_u^n}{n}\frac{S_u^n}{n}\\
-\lambda(E\rho)\frac{I_u^n}{n}\frac{S_u^n}{n} & \lambda(E\rho)\frac{I_u^n}{n}\frac{S_u^n}{n}+\frac{I_u^n}{n}
\end{pmatrix},
\]
and
\[
b_u^n=\sum_{i=1}^2l_i(\nabla^TH_i)(\xi_u^n)=
\begin{pmatrix}
-\lambda(E\rho)i_u(\xi^n) & -\lambda(E\rho)s_u(\xi^n)\\
\lambda(E\rho)i_u(\xi^n) & \lambda(E\rho)s_u(\xi^n)-1
\end{pmatrix}
\]
while $\xi^n_u$ is a convex combination of $\vartheta_u^n$ and $\widehat{x}_u$ satisfying $H_2(\vartheta_u^n)-H_2(\widehat{x}_u)=\nabla H_2(\xi^n_u)\cdot(\vartheta_u^n-\widehat{x}_u)$, the existence of which follows from Lagrange's mean value theorem.

By Lemma \ref{lemma 5.2 XiIsMartingale}, we define $\widehat{Q}_{\lambda,n}^{\omega,g}$ as the quenched measure such that
\[
\frac{d\widehat{Q}_{\lambda,n}^{\omega,g}}{dP_{\lambda,n}^\omega}=\Xi_{T_0}^n(g)
\]
for any $\omega\in X$ and $g\in C^2\big([0,T_0], \mathbb{R}^2\big)$. We further define $\widehat{Q}^g_{\lambda,n}$ as the annealed measure such that
\[
\widehat{Q}^g_{\lambda,n}(\cdot)=\int_{X}\widehat{Q}_{\lambda,n}^{\omega,g}(\cdot)P(d\omega).
\]
For $x,y\geq 0$ such that $x+y\leq 1$, we define
\[
\widehat{Q}^{g,x,y}_{\lambda,n}(\cdot)=\widehat{Q}^g_{\lambda,n}\Big(\cdot\Big|S_0^n=\lfloor np_0+a_nx\rfloor, I_0^n=\lfloor np_1+a_ny\rfloor\Big).
\]
Then, we have the following lemma, which is an analogue of Lemma \ref{lemma 3.3 LLNunderChangeMeasure} and crucial for the proof of Equation \eqref{equ MDP lower bound}.

\begin{lemma}\label{lemma 5.3 MDPLLNunderChangeMeasure}
For any $x,y\geq 0$ such that $x+y\leq 1$ and any $g\in C^2\big([0,T_0], \mathbb{R}^2\big)$, $\nu^n$
converges in $\widehat{Q}^{g,x,y}_{\lambda,n}$-probability to the solution
\[\{\overline{x}_t=(\overline{s}_t,\overline{i}_t)^T:~0\leq t\leq T_0\}\]
to the ODE
\[
\begin{cases}
&\frac{d}{dt}\overline{x}_t=b_t\overline{x}_t+\sigma_tg_t, \\
&(\overline{s}_0,\overline{i}_0)=(x,y).
\end{cases}
\]
\end{lemma}
The following proof of Lemma \ref{lemma 5.3 MDPLLNunderChangeMeasure} is similar with that of Lemma \ref{lemma 3.3 LLNunderChangeMeasure}, where the generalized version of Girsanov's theorem is utilized.

\proof[Outline of the proof of Lemma \ref{lemma 5.3 MDPLLNunderChangeMeasure}]

Since $S_t^n, I_t^n\leq n$, there exists $K_8$ not depending on $n$ such that $\Xi_t^n(g)\leq e^{\frac{a_n^2}{n}K_8}$ for all $0\leq t\leq T_0$, then according to the analysis in the proofs of Lemmas \ref{lemma 2.3 LLNofSIR}, \ref{lemma 3.1}, \ref{lemma 5.1} and Cauchy-Schwartzs's inequality,
\[
\sup_{0\leq u\leq T_0}|\vartheta_u^n-\widehat{x}_u|, \sup_{0\leq u\leq T_0}|\varepsilon_u^n|, \sup_{0\leq u\leq T_0}|b_u^n-b_u| \text{~and~}
\sup_{0\leq u\leq T_0}|\sigma^n_u-\sigma_u|
\]
converges to $0$ in both $P_{\lambda,n}$-probability and $\widehat{Q}^{g,x,y}_{\lambda,n}$-probability as $n\rightarrow+\infty$.

Then, by Equation \eqref{equ 5.1 ExpressionXi},
\begin{align}\label{equ 5.2 ExpressionXiSecondversion}
\Xi_t^n(g)=\exp\Bigg\{&\frac{a_n^2}{n}\Bigg(g_t\cdot \nu_t^n-g_0\cdot (x,y)^T-\int_0^tg_u^\prime\cdot \nu_u^n+(b_u\nu_n^n)\cdot g_u^n \notag\\
&+\frac{1}{2}g_u^T\sigma_u g_udu+o(1)\Bigg)\Bigg\}
\end{align}
under both probability measures.

Similar with that in  the proof of Lemma \ref{lemma 3.3 LLNunderChangeMeasure}, we define
\[
\zeta_1(t,\eta)=\sum_{i=1}^n1_{\{\eta(i)=0\}}-n\widehat{s}_t \text{~and~}\zeta_2(t,\eta)=\sum_{i=1}^n1_{\{\eta(i)=1\}}-n\widehat{i}_t.
\]
Hence, $\zeta_1(t,\eta_t^n)=S_t^n-n\widehat{s}_t$ while $\zeta_2(t,\eta_t^n)=I_t^n-n\widehat{i}_t$. We further define
\[
\mathcal{M}_t(V_g)=V_g(t,\eta_t^n)-V_g(0,\eta_0^n)-\int_0^t (\frac{\partial}{\partial u}+\Omega_n)V_g(u, \eta^n_u)du,
\]
\[
\widehat{\mathcal{M}}_t(V_g)=\int_0^t \frac{1}{V_g(u-,\eta_{u-}^n)}d\mathcal{M}_{u}(V_g)
\]
for $g\in C^2\big([0,T_0], \mathbb{R}^2\big)$ and
\[
\mathcal{M}_t(\zeta_i)=\zeta_i(t,\eta_t^n)-\zeta_i(0,\eta_0^n)-\int_0^t (\frac{\partial}{\partial u}+\Omega_n)\zeta_i(u,\eta_u^n)du,
\]
\[
\widetilde{\mathcal{M}}_t(\zeta_i)=\mathcal{M}_t(\zeta_i)-<\mathcal{M}(\zeta_i), \widehat{\mathcal{M}}(V_g)>_t
\]
for $i=1,2$.

Then, according to the generalized version of Girsanov's theorem and a similar analysis with that in the proof of Lemma \ref{lemma 3.3 LLNunderChangeMeasure}, we have following analogue of Equations \eqref{equ StLDP} and \eqref{equ ItLDP},
\begin{equation}\label{equ ItStMDP}
\begin{cases}
&\nu_t^n=\nu_0^n+\int_0^t (b_u+o(1))\nu_u^n+(\sigma_u+o(1))g_udu+\frac{1}{a_n}(\widetilde{\mathcal{M}}_t(\zeta_1),\widetilde{\mathcal{M}}_t(\zeta_2))^T, \\
&\nu_0^n=\big(\frac{\lfloor np_0+a_nx\rfloor-np_0}{a_n},\frac{\lfloor np_1+a_ny\rfloor-np_1}{a_n}\big)^T.
\end{cases}
\end{equation}
Note that, to obtain Equation \eqref{equ ItStMDP}, we should utilize Equation \eqref{equ 3.1} and the fact that $\varepsilon_u^n, \sigma_u^n-\sigma_u$ converges to $0$ to check that
\[
\big(<\mathcal{\mathcal{M}}(\zeta_1), \widehat{\mathcal{M}}(V_g)>_t,~<\mathcal{M}(\zeta_2), \widehat{M}(V_g)>_t\big)^T=(\sigma_t+o(1))g_t.
\]
Since the calculation is not difficult but a little tedious, we omit details here.

With Equation \eqref{equ ItStMDP}, we only need to show that $\frac{1}{a_n}\sup_{0\leq t\leq T_0}|\widetilde{\mathcal{M}}_t(\zeta_i)|$ converges to $0$ in $\widehat{Q}^{g,x,y}_{\lambda,n}$-probability as $n\rightarrow+\infty$ to complete this proof. To check this property, we only need to show that
\[
\frac{1}{a_n^2}\big[\widetilde{\mathcal{M}}(\zeta_i),\widetilde{\mathcal{M}}(\zeta_i)\big]_{T_0}=
\frac{1}{a_n^2}\big[\mathcal{M}(\zeta_i),\mathcal{M}(\zeta_i)\big]_{T_0}
 \]
converges to $0$ in $\widehat{Q}^{g,x,y}_{\lambda,n}$-probability, which holds according to a similar analysis with that at the end of the proof of Lemma \ref{lemma 3.3 LLNunderChangeMeasure}. In detail, since $S_t^n, I_t^n\leq n$, $\big[\mathcal{M}(\zeta_i),\mathcal{M}(\zeta_i)\big]_{T_0}$ is stochastically dominated from above by $K_9Y_{nK_{10}T_0}$ under $P_{\lambda,n}$, where $K_9,K_{10}\in (0,+\infty)$ does not depend on $n$ and $\{Y_t\}_{t\geq 0}$ is the Poisson process with rate one. Therefore, by chebyshev's inequality and the fact that $\Xi_{T_0}^n(g)\leq e^{K_8\frac{a_n^2}{n}}$,
\[
\widehat{Q}^{g,x,y}_{\lambda,n}\Big(\frac{1}{a_n^2}\big[\widetilde{\mathcal{M}}(\zeta_i),\widetilde{\mathcal{M}}(\zeta_i)\big]_{T_0}\geq \epsilon\Big)\leq e^{-a_n^2\epsilon}e^{K_{8}\frac{a_n^2}{n}}e^{nK_{10}T_0(e^{K_9}-1)}\rightarrow 0
\]
for any $\epsilon>0$ according to $\frac{a_n^2}{n}\rightarrow +\infty$ and the proof is complete.

\qed

Now we give the proof of Equation \eqref{equ MDP lower bound}, which is similar with that of Equation \eqref{equ LDP lower}.

\proof[Outline of the proof of Equation \eqref{equ MDP lower bound}]

We only need to deal with the case where $\inf_{f\in O}(J_{dyn}(f)+J_{ini}(f_0))<+\infty$. For any $\epsilon>0$, let $\widetilde{f}^\epsilon\in O$ such that
\[
J_{dyn}(\widetilde{f}^\epsilon)+J_{ini}(\widetilde{f}^\epsilon_0)<\inf_{f\in O}(J_{dyn}(f)+J_{ini}(f_0))+\epsilon.
\]
Then, according to the analysis in the proof of Theorem \ref{theorem 2.5 JdynEpress}, there exists $\widetilde{g}^\epsilon$ such that
\[
(\widetilde{f}^\epsilon)^\prime_t=b_t\widetilde{f}^\epsilon_t+\sigma_t\widetilde{g}^\epsilon_t \text{~and~}J_{dyn}(\widetilde{f}^\epsilon)
=\frac{1}{2}\int_0^{T_0}(\widetilde{g}^\epsilon)^T_t\sigma_t\widetilde{g}^\epsilon_tdt.
\]
For $n\geq 1$, let $\widetilde{g}^n\in C^2\big([0,T_0], \mathbb{R}^2\big)$ such that $\widetilde{g}^n$ converges to $\widetilde{g}^\epsilon$ in $L^2$-norm and let $\widetilde{f}^n$ be the solution to the ode
\[
\begin{cases}
&\frac{d}{dt}\widetilde{f}^n_t=b_t\widetilde{f}^n_t+\sigma_t\widetilde{g}^n_t,\\
&\widetilde{f}^n_0=\widetilde{f}^\epsilon_0.
\end{cases}
\]
Then, $\widetilde{f}^n$ converges to $\widetilde{f}^\epsilon$ uniformly on $[0,T_0]$ and
\[
J_{dyn}(\widetilde{f}^n)=\frac{1}{2}\int_0^{T_0}(\widetilde{g}^n)^T_t\sigma_t\widetilde{g}^n_tdt
\rightarrow \frac{1}{2}\int_0^{T_0}(\widetilde{g}^\epsilon)^T_t\sigma_t\widetilde{g}^\epsilon_tdt=J_{dyn}(\widetilde{f}^\epsilon).
\]
Hence, there exists integer $m_1>1$ such that $\widetilde{f}^{m_1}\in O$ and $|J_{dyn}(\widetilde{f}^{m_1})-J_{dyn}(\widetilde{f}^\epsilon)|<\epsilon$.

According to Equation \eqref{equ 5.1 ExpressionXi} and the fact that
\begin{align*}
J_{dyn}(\widetilde{f}^{m_1})=&\frac{1}{2}\int_0^{T_0}(\widetilde{g}^{m_1})^T_t\sigma_t\widetilde{g}^{m_1}_tdt\\
=&\widetilde{f}^{m_1}_{T_0}\cdot \widetilde{g}^{m_1}_{T_0}-\widetilde{f}^{m_1}_0\cdot \widetilde{g}^{m_1}_0-\int_0^{T_0}\widetilde{f}^{m_1}_t\cdot (\widetilde{g}^{m_1}_t)^\prime dt\\
&-\int_0^{T_0}(b_t\widetilde{f}^{m_1}_t)\cdot \widetilde{g}^{m_1}_tdt-\frac{1}{2}\int_0^{T_0}(\widetilde{g}^{m_1}_t)^T\sigma_t \widetilde{g}^{m_1}_tdt,
\end{align*}
there exists $\widetilde{\delta}>0$ such that $B(\widetilde{f}^{m_1},\widetilde{\delta})\subseteq O$ and
\[
\Xi_{T_0}^n({\widetilde{g}^{m_1}})\leq \exp\Big\{\frac{a_n^2}{n}\big(J_{dyn}(\widetilde{f}^{m_1})+\epsilon\big)\Big\}
\]
conditioned on $\nu^n\in B(\widetilde{f}^{m_1},\widetilde{\delta})$, $\frac{n\delta_n}{a_n}\leq \widetilde{\delta}$ and $\|\vartheta^n-\widehat{x}\|\leq \widetilde{\delta}$. As a result,
\begin{align*}
P_{\lambda,n}\big(\nu^n\in O\big)
\geq e^{-\frac{a_n^2}{n}(J_{dyn}(\widetilde{f}^{m_1})+\epsilon)}
\widehat{Q}^{\widetilde{g}^{m_1}}_{\lambda,n}\Big(\nu^n\in B(\widetilde{f}^{m_1},\widetilde{\delta}), \frac{n\delta_n}{a_n}\leq \widetilde{\delta},
\|\vartheta^n-\widehat{x}\|\leq\widetilde{\delta}\Big).
\end{align*}
According to the analysis in the proofs of Lemmas \ref{lemma 2.3 LLNofSIR}, \ref{lemma 5.1} and the fact that $\Xi^n_{T_0}(\widetilde{g}^{m_1})\leq e^{\frac{a_n^2}{n}K_8}$,
\[
\lim_{n\rightarrow+\infty}\frac{n}{a_n^2}\log \widehat{Q}^{\widetilde{g}^{m_1}}_{\lambda,n}\big(\frac{n\delta_n}{a_n}>\widetilde{\delta}\text{~or~}
\|\vartheta^n-\widehat{x}\|>\widetilde{\delta}\big)=-\infty.
\]
Hence,
\begin{align*}
&\liminf_{n\rightarrow+\infty}\frac{n}{a_n^2}\log \widehat{Q}^{\widetilde{g}^{m_1}}_{\lambda,n}\big(\frac{n\delta_n}{a_n}\leq \widetilde{\delta},
\|\vartheta^n-\widehat{x}\|\leq \widetilde{\delta}, \nu^n\in B(\widetilde{f}^{m_1},\widetilde{\delta})\big) \\
&=\liminf_{n\rightarrow+\infty}\frac{n}{a_n^2}\log \widehat{Q}^{\widetilde{g}^{m_1}}_{\lambda,n}\big(\nu^n\in B(\widetilde{f}^{m_1},\widetilde{\delta})\big)
\end{align*}
and
\begin{equation}\label{equ 5.4}
\liminf_{n\rightarrow+\infty}\frac{n}{a_n^2}\log P_{\lambda,n}\big(\nu^n\in O\big)
\geq -J_{dyn}(\widetilde{f}^{m_1})-\epsilon+\liminf_{n\rightarrow+\infty}\frac{n}{a_n^2}\log \widehat{Q}^{\widetilde{g}^{m_1}}_{\lambda,n}\big(\nu^n\in B(\widetilde{f}^{m_1},\widetilde{\delta})\big).
\end{equation}
By Lemma \ref{lemma 5.3 MDPLLNunderChangeMeasure}, since $(\widetilde{f}^{m_1})^\prime_t=b_t\widetilde{f}^{m_1}_t+\sigma_t\widetilde{g}^{m_1}_t$,
\[
\lim_{n\rightarrow+\infty}\widehat{Q}^{\widetilde{g}^{m_1},s_0(\widetilde{f}^\epsilon),i_0(\widetilde{f}^\epsilon)}_{\lambda,n}\big(\nu^n\in B(\widetilde{f}^{m_1},\widetilde{\delta})\big)=1.
\]
As a result,
\begin{align*}
&\liminf_{n\rightarrow+\infty}\frac{n}{a_n^2}\log \widehat{Q}^{\widetilde{g}^{m_1}}_{\lambda,n}\big(\nu^n\in B(\widetilde{f}^{m_1},\widetilde{\delta})\big) \\
&\geq \liminf_{n\rightarrow+\infty}\frac{n}{a_n^2}\log \widehat{Q}^{\widetilde{g}^{m_1}}_{\lambda,n}\big(\nu^n\in B(\widetilde{f}^{m_1},\widetilde{\delta}), S_0^n=\lfloor np_0+a_ns_0(\widetilde{f}^\epsilon)\rfloor, I_0^n=\lfloor np_1+a_ni_0(\widetilde{f}^\epsilon)\rfloor\big)\\
&=\liminf_{n\rightarrow+\infty}\frac{n}{a_n^2}\log \Big[\widehat{Q}^{\widetilde{g}^{m_1},s_0(\widetilde{f}^\epsilon),i_0(\widetilde{f}^\epsilon)}_{\lambda,n}\big(\nu^n\in B(\widetilde{f}^{m_1},\widetilde{\delta})\big)\\
&\text{~\quad\quad}\times P_{\lambda,n}\big(S_0^n=\lfloor np_0+a_ns_0(\widetilde{f}^\epsilon)\rfloor, I_0^n=\lfloor np_1+a_ni_0(\widetilde{f}^\epsilon)\rfloor\big)\Big]\\
&=\liminf_{n\rightarrow+\infty}\frac{n}{a_n^2}\log P_{\lambda,n}\big(S_0^n=\lfloor np_0+a_ns_0(\widetilde{f}^\epsilon)\rfloor, I_0^n=\lfloor np_1+a_ni_0(\widetilde{f}^\epsilon)\rfloor\big).
\end{align*}
For given $x,y\in \mathbb{R}$, according to Strling's formula and the definition of $M_{_0}$,
\begin{align*}
&\liminf_{n\rightarrow+\infty}\frac{n}{a_n^2}\log P_{\lambda,n}\big(S_0^n=\lfloor np_0+a_nx\rfloor, I_0^n=\lfloor np_1+a_ny\rfloor\big)\\
&=\liminf_{n\rightarrow+\infty}\frac{n}{a_n^2}\log \Big[{n\choose \lfloor np_0+a_nx\rfloor}
{n-\lfloor np_0+a_nx\rfloor \choose \lfloor np_1+a_ny\rfloor}\\
&\text{\quad\quad}\times p_0^{\lfloor np_0+a_nx\rfloor}p_1^{\lfloor np_1+a_ny\rfloor}(1-p_0-p_1)^{n-\lfloor np_0+a_nx\rfloor-\lfloor np_1+a_ny\rfloor}\Big]\\
&=-\frac{1}{2}\big[\frac{1}{p_0}x^2+\frac{1}{p_1}y^2+\frac{(x+y)^2}{1-p_0-p_1}\big]=-\frac{1}{2}(x,y)M^{-1}_{_0}(x,y)^T.
\end{align*}
Therefore,
\[
\liminf_{n\rightarrow+\infty}\frac{n}{a_n^2}\log \widehat{Q}^{\widetilde{g}^{m_1}}_{\lambda,n}\big(\nu^n\in B(\widetilde{f}^{m_1},\widetilde{\delta})\big)\geq -\frac{1}{2}(\widetilde{f}^\epsilon_0)^TM^{-1}_{_0}\widetilde{f}^\epsilon_0
=-J_{ini}(\widetilde{f}^\epsilon_0)
\]
by Theorem \ref{theorem 2.5 JdynEpress}. Then, by Equation \eqref{equ 5.4},
\begin{align*}
\liminf_{n\rightarrow+\infty}\frac{n}{a_n^2}\log P_{\lambda,n}\big(\nu^n\in O\big)
&\geq -J_{dyn}(\widetilde{f}^{m_1})-\epsilon-J_{ini}(\widetilde{f}^\epsilon_0)\\
&\geq -\big(J_{dyn}(\widetilde{f}^\epsilon)+J_{ini}(\widetilde{f}^\epsilon_0)\big)-2\epsilon \\
&\geq -\inf_{f\in O}(J_{dyn}(f)+J_{ini}(f_0))-3\epsilon,
\end{align*}
Equation \eqref{equ MDP lower bound} follows from which directly since $\epsilon$ is arbitrary.

\qed

The proof of Equation \eqref{equ MDP upper bound} is similar with that of Equation \eqref{equ LDP upper}. First we show that the the Equation holds for compact sets.

\begin{lemma}\label{lemma 5.4}
For any closed set $C\subseteq \mathcal{D}\big([0,T_0], \mathbb{R}^2\big)$,
\[
\limsup_{n\rightarrow+\infty}\frac{n}{a_n^2}\log P_{\lambda,n}\big(\nu^n\in C\big)\leq-\inf_{f\in C}(J_{dyn}(f)+J_{ini}(f_0)).
\]
\end{lemma}

\proof[Proof of Lemma \ref{lemma 5.4}]

For any $f\in \mathcal{D}\big([0,T_0], \mathbb{R}^2\big)$ and $g\in C^2\big([0,T_0], \mathbb{R}^2\big)$, let
\[
\mathcal{L}_{f}(g)=f_{T_0}\cdot g_{T_0}-f_0\cdot g_0-\int_0^{T_0}f_t\cdot g_t^\prime dt-\int_0^{T_0}(b_tf_t)\cdot g_tdt-\frac{1}{2}\int_0^{T_0}g_t^T\sigma_t g_tdt,
\]
then $J_{dyn}(f)=\sup_{g\in C^2\big([0,T_0], \mathbb{R}^2\big)}\mathcal{L}_f(g)$.

By Equation \eqref{equ 5.1 ExpressionXi}, for given $\epsilon>0$, $g\in C^2\big([0,T_0], \mathbb{R}^2\big)$ and compact set $C\subseteq \mathcal{D}\big([0,T_0], \mathbb{R}^2\big)$, there exists $\widetilde{\delta}_2$ depending on $\epsilon, g, C$ such that
\[
\Xi^n_{T_0}(g)\geq \exp\Big\{\frac{a_n^2}{n}\big(\mathcal{L}_{\nu^n}(g)-\epsilon\big)\Big\}
\]
for sufficiently large $n$ conditioned on $\frac{n\delta_n}{a_n}\leq \widetilde{\delta}_2$, $\|\vartheta^n-\widehat{x}\|\leq \widetilde{\delta}_2$
and $\nu^n\in C$. Therefore, for any $y\in \mathbb{R}^2$, according to the fact that $\{\Xi_t^n(g)\}_{0\leq t\leq T_0}$ is a martingale,
\begin{align}\label{equ 5.5}
&E_{\lambda,n}e^{\frac{a^2_n}{n}y\cdot \nu_0^n}=E_{\lambda,n}\big(e^{\frac{a^2_n}{n}y\cdot \nu_0^n}\Xi_{T_0}^n(g)\big)\notag \\
&\geq E\big(e^{\frac{a^2_n}{n}y\cdot \nu_0^n}\Xi_{T_0}^n(g)1_{\{\nu^n\in C, \|\vartheta^n-\widehat{x}\|\leq \widetilde{\delta}_2, \frac{n\delta_n}{a_n}\leq \widetilde{\delta}_2\}}\big) \\
&\geq \exp\Big\{\frac{a^2_n}{n}\big(\inf_{f\in C}\{y\cdot f_0^n+\mathcal{L}_{f}(g)\}-\epsilon\big)\Big\}P_{\lambda,n}\big(f\in C, \|\vartheta^n-\widehat{x}\|\leq \widetilde{\delta}_2, \frac{n\delta_n}{a_n}\leq \widetilde{\delta}_2\big). \notag
\end{align}
By Assumption (\textbf{A}) and Taylor's expansion formula,
\begin{align*}
E_{\lambda,n}e^{\frac{a^2_n}{n}y\cdot \nu_0^n}&=e^{-\frac{a_n(s_yp_0+i_yp_1)}{n}}\big(1-p_0-p_1+p_1e^{\frac{a_n}{n}i_y}+p_0e^{\frac{a_n}{n}s_y}\big)^n \\
&=\exp\big\{\frac{a_n^2}{2n}(y^TM_{_0}y+o(1))\big\}.
\end{align*}
Consequently, by Equation \eqref{equ 5.5},
\begin{align*}
&\limsup_{n\rightarrow+\infty}\frac{n}{a_n^2}\log P_{\lambda,n}\big(\nu^n\in C, \|\vartheta^n-\widehat{x}\|\leq \widetilde{\delta}_2, \frac{n\delta_n}{a_n}\leq \widetilde{\delta}_2\big)\\
&\leq -\inf_{f\in C}\{y\cdot f_0-\frac{1}{2}y^TM_{_0}y+\mathcal{L}_{f}(g)\}+\epsilon.
\end{align*}
According to the analysis in the proof of Lemmas \ref{lemma 2.3 LLNofSIR} and \ref{lemma 5.1},
\begin{align*}
\limsup_{n\rightarrow+\infty}\frac{n}{a_n^2}\log P_{\lambda,n}\big(\nu^n\in C, \|\vartheta^n-\widehat{x}\|\leq \widetilde{\delta}_2, \frac{n\delta_n}{a_n}\leq \widetilde{\delta}_2\big)=\limsup_{n\rightarrow+\infty}\frac{n}{a_n^2}\log P_{\lambda,n}\big(\nu^n\in C\big)
\end{align*}
and hence
\[
\limsup_{n\rightarrow+\infty}\frac{n}{a_n^2}\log P_{\lambda,n}\big(\nu^n\in C\big)\\
\leq -\inf_{f\in C}\{y\cdot f_0-\frac{1}{2}y^TM_{_0}y+\mathcal{L}_{f}(g)\}+\epsilon.
\]
Since $\epsilon$ is arbitrary,
\[
\limsup_{n\rightarrow+\infty}\frac{n}{a_n^2}\log P_{\lambda,n}\big(\nu^n\in C\big)\\
\leq -\inf_{f\in C}\{y\cdot f_0-\frac{1}{2}y^TM_{_0}y+\mathcal{L}_{f}(g)\}.
\]
Since $y, g$ are arbitrary,
\begin{equation}\label{equ 5.6}
\limsup_{n\rightarrow+\infty}\frac{n}{a_n^2}\log P_{\lambda,n}\big(\nu^n\in C\big)\\
\leq -\sup_{y,g}\inf_{f\in C}\{y\cdot f_0-\frac{1}{2}y^TM_{_0}y+\mathcal{L}_{f}(g)\}.
\end{equation}
Since $C$ is compact, according to the minimax theorem,
\begin{align*}
&\sup_{y,g}\inf_{f\in C}\{y\cdot f_0-\frac{1}{2}y^TM_{_0}y+\mathcal{L}_{f}(g)\}
=\inf_{f\in C}\sup_{y,g}\{y\cdot f_0-\frac{1}{2}y^TM_{_0}y+\mathcal{L}_{f}(g)\} \\
&=\inf_{f\in C}\big(\sup_{y}\{y\cdot f_0-\frac{1}{2}y^TM_{_0}y\}+\sup_{g}\mathcal{L}_{f}(g)\big)
=\inf_{f\in C}\big(J_{ini}(f_0)+J_{dyn}(f)\big),
\end{align*}
Lemma \ref{lemma 5.4} follows from which and Equation \eqref{equ 5.6} directly.

\qed

At last, we give the proof of Equation \eqref{equ MDP upper bound}.

\proof[Outline of the proof of Equation \eqref{equ MDP upper bound}]

With Lemma \ref{lemma 5.4}, we only need to show that $\{\nu^n\}_{n\geq 1}$ are exponential tight, which is equivalent to the following two properties.

(\textbf{1}) \[
\limsup_{M\rightarrow+\infty}\limsup_{n\rightarrow+\infty}\frac{n}{a_n^2}\log P_{\lambda,n}\big(\|\nu^n\|>M\big)=-\infty.
\]
(\textbf{2}) For any $\epsilon>0$,
\[
\limsup_{\delta\rightarrow 0}\limsup_{n\rightarrow+\infty}\frac{n}{a_n^2}\log \sup_{\tau \in \mathcal{T}_0}
P_{\lambda,n}\big(\sup_{0\leq t\leq \delta}|\nu^n_{t+\tau}-\nu^n_\tau|_1>\epsilon\big)=-\infty,
\]
where $\mathcal{T}_0$ is the set of stopping times of $\{\eta_t^n\}_{0\leq t\leq T_0}$ with upper bound $T_0$ and $|x|_1$ is the $l_1$-norm of $x\in \mathbb{R}^2$.

By Lemmas \ref{lemma 3.1} and \ref{lemma 5.1}, Properties (\textbf{1}) and (\textbf{2}) are equivalent to

(\textbf{\.{1}}) \[
\limsup_{M\rightarrow+\infty}\limsup_{n\rightarrow+\infty}\frac{n}{a_n^2}\log P_{\lambda,n}\big(\|\nu^n\|>M,\delta_n\leq 1\big)=-\infty.
\]
(\textbf{\.{2}}) For any $\epsilon>0$,
\[
\limsup_{\delta\rightarrow 0}\limsup_{n\rightarrow+\infty}\frac{n}{a_n^2}\log \sup_{\tau \in \mathcal{T}_0}
P_{\lambda,n}\big(\sup_{0\leq t\leq \delta}|\nu^n_{t+\tau}-\nu^n_\tau|_1>\epsilon, \frac{n\delta_n }{a_n}\leq 1\big)=-\infty.
\]

To check (\textbf{\.{1}}), we utilize the analysis introduced in Chapter 11 of \cite{Ethier1986}. Since this is a well-known analysis, we only give an outline. According to the the generator $\Omega_n$ of our process, $S_t^n-n\widehat{s}_t$ and $I_t^n-n\widehat{i}_t$ can be written as
\begin{align*}
S_t^n-n\widehat{s}_t=&S_0^n-n\widehat{s}_0-\overline{Y}_2\big(\int_0^tn\lambda(E\rho\frac{S_u^n}{n}\frac{I_u^n}{n}+\varepsilon_u^n) du\big)\\
&-\int_0^t n\lambda E\rho(\frac{S^n_u}{n}\frac{I^n_u}{n}-\widehat{s}_u\widehat{i}_u) du-\int_0^tn\lambda \varepsilon_u^n du
\end{align*}
and
\begin{align*}
I_t^n&-n\widehat{i}_t=I_0^n-n\widehat{i}_0-\overline{Y}_1\big(\int_0^t I_u^ndu\big)-\int_0^t n(\frac{I_u^n}{n}-\widehat{i}_u)du \\
&+\overline{Y}_2\big(\int_0^tn\lambda(E\rho\frac{S_u^n}{n}\frac{I_u^n}{n}+\varepsilon_u^n) du\big)+\int_0^tn\lambda E\rho(\frac{S_u^n}{n}\frac{I_u^n}{n}-\widehat{s}_u\widehat{i}_u)du+\int_0^tn\varepsilon_u^n\lambda du,
\end{align*}
where $\varepsilon_u^n=\frac{\gamma(\mathcal{S}_u^n, \mathcal{I}_u^n)-(E\rho)S_u^nI_u^n}{n^2}$ and $\overline{Y}_i(t)=Y_i(t)-t$ such that $\{Y_i(t)\}_{t\geq 0}$ is a Poisson process with rate one for $i=1,2$. Then, according to Grownwall's inequality and the fact that $S_t^n, I_t^n\leq n$,
\[
\|\nu^n_t\|_1\leq \varepsilon_0e^{K_{11}T_0}
\]
for all $0\leq t\leq T_0$ conditioned on $\delta_n\leq 1$, where $K_{11}\in (0,+\infty)$ does not depend on $n$ due to the Lipschitz's condition of $H_1, H_2$ defined in Section \ref{section two} while
\begin{align*}
\varepsilon_0=& \frac{1}{a_n}\Big(\sup_{0\leq t\leq nT_0}|\overline{Y}_1(t)|+2\sup_{0\leq t\leq nK_{12}T_0}|\overline{Y}_2(t)|\Big) \\
&+\frac{1}{a_n}\big(|I_0^n-n\widehat{i}_0|+|S_0^n-n\widehat{s}_0|\big)+\frac{n\lambda T_0}{a_n}\delta_n,
\end{align*}
where $K_{12}=\lambda(E\rho+1)$. Consequently, Property (\textbf{\.{1}}) holds according to well known moderate deviation principles of Poisson processes and sums of i.i.d. random variables and Lemma \ref{lemma 5.1}.

Now we only need to check Property (\textbf{\.{2}}). With Property (\textbf{\.{1}}), we only need to check that
\begin{align}\label{equ 5.7}
\limsup_{\delta\rightarrow 0}\limsup_{n\rightarrow+\infty}\frac{n}{a_n^2}\log \sup_{\tau \in \mathcal{T}_0}
P_{\lambda,n}\big(&\sup_{0\leq t\leq \delta}|\frac{S^n_{t+\tau}-n\widehat{s}_{t+\tau}}{a_n}-\frac{S^n_\tau-n\widehat{s}_\tau}{a_n}|>\epsilon,\notag\\ &\frac{n\delta_n }{a_n}\leq 1, \|\nu^n\|\leq M\big)=-\infty
\end{align}
and
\begin{align}\label{equ 5.8}
\limsup_{\delta\rightarrow 0}\limsup_{n\rightarrow+\infty}\frac{n}{a_n^2}\log \sup_{\tau \in \mathcal{T}_0}
P_{\lambda,n}\big(&\sup_{0\leq t\leq \delta}|\frac{I^n_{t+\tau}-n\widehat{i}_{t+\tau}}{a_n}-\frac{I^n_\tau-n\widehat{i}_\tau}{a_n}|>\epsilon,\notag\\ &\frac{n\delta_n }{a_n}\leq 1, \|\nu^n\|\leq M\big)=-\infty
\end{align}
for any $M>0$ and $\epsilon>0$.

Let $e_1=(1,0)^T$ and $\chi_1(t)\equiv e_1$, then for any $\theta>0$ and $\tau \in \mathcal{T}_0$, $\{\frac{\Xi^n_{\tau+t}(\theta\chi_1)}{\Xi^n_\tau(\theta\chi_1)}\}_{t\geq 0}$ is a martingale according to Lemma \ref{lemma 5.2 XiIsMartingale}. By Equation \eqref{equ 5.1 ExpressionXi}, conditioned on $\frac{n\delta_n }{a_n}\leq 1$ and $\|\nu^n\|\leq M$, there exists $K_{13}, K_{14}\in (0,+\infty)$ not depending on $n$ such that
\begin{align*}
\frac{\Xi^n_{\tau+t}(\theta\chi_1)}{\Xi^n_\tau(\theta\chi_1)}
\geq \exp\Bigg\{\frac{a_n^2}{n}\Big[\theta\big(\frac{S^n_{t+\tau}-n\widehat{s}_{t+\tau}}{a_n}-\frac{S^n_\tau-n\widehat{s}_\tau}{a_n}\big)                           -\delta(\theta K_{13}+\theta^2K_{14})\Big]\Bigg\}
\end{align*}
for sufficiently large $n$ and any $0\leq t\leq \delta$. Then,
\begin{align*}
&\Big\{\sup_{0\leq t\leq \delta}\big(\frac{S^n_{t+\tau}-n\widehat{s}_{t+\tau}}{a_n}-\frac{S^n_\tau-n\widehat{s}_\tau}{a_n}\big)>\epsilon, \frac{n\delta_n }{a_n}\leq 1, \|\nu^n\|\leq M\Big\} \\
& \subseteq \Big\{\sup_{0\leq t\leq \delta}\frac{\Xi^n_{\tau+t}(\theta\chi_1)}{\Xi^n_\tau(\theta\chi_1)}
\geq e^{\frac{a_n^2}{n}[\theta\epsilon-\delta(\theta K_{13}+\theta^2K_{14})]}\Big\}.
\end{align*}
By Doob's inequality,
\begin{align*}
P_{\lambda,n}\big(\sup_{0\leq t\leq \delta}\frac{\Xi^n_{\tau+t}(\theta\chi_1)}{\Xi^n_\tau(\theta\chi_1)}
\geq e^{\frac{a_n^2}{n}[\theta\epsilon-\delta(\theta K_{13}+\theta^2K_{14})]}\big)
&\leq e^{-\frac{a_n^2}{n}[\theta\epsilon-\delta(\theta K_{13}+\theta^2K_{14})]}E_{\lambda,n}\frac{\Xi^n_{\tau+t}(\theta\chi_1)}{\Xi^n_\tau(\theta\chi_1)}\\
&=e^{-\frac{a_n^2}{n}[\theta\epsilon-\delta(\theta K_{13}+\theta^2K_{14})]}
\end{align*}
and hence
\begin{align*}
\limsup_{\delta\rightarrow 0}\limsup_{n\rightarrow+\infty}\frac{n}{a_n^2}\log \sup_{\tau \in \mathcal{T}_0}
P_{\lambda,n}\big(&\sup_{0\leq t\leq \delta}(\frac{S^n_{t+\tau}-n\widehat{s}_{t+\tau}}{a_n}-\frac{S^n_\tau-n\widehat{s}_\tau}{a_n})>\epsilon,\notag\\ &\frac{n\delta_n }{a_n}\leq 1, \|\nu^n\|\leq M\big)\leq -\theta\epsilon.
\end{align*}
Since $\theta$ is arbitrary,
\begin{align}\label{equ 5.9}
\limsup_{\delta\rightarrow 0}\limsup_{n\rightarrow+\infty}\frac{n}{a_n^2}\log \sup_{\tau \in \mathcal{T}_0}
P_{\lambda,n}\big(&\sup_{0\leq t\leq \delta}(\frac{S^n_{t+\tau}-n\widehat{s}_{t+\tau}}{a_n}-\frac{S^n_\tau-n\widehat{s}_\tau}{a_n})>\epsilon,\notag\\ &\frac{n\delta_n }{a_n}\leq 1, \|\nu^n\|\leq M\big)=-\infty.
\end{align}
Since $\{\frac{\Xi^n_{\tau+t}(-\theta\chi_1)}{\Xi^n_\tau(-\theta\chi_1)}\}_{t\geq 0}$ is also a martingale for any $\theta>0$, according to a similar analysis,
\begin{align}\label{equ 5.10}
\limsup_{\delta\rightarrow 0}\limsup_{n\rightarrow+\infty}\frac{n}{a_n^2}\log \sup_{\tau \in \mathcal{T}_0}
P_{\lambda,n}\big(&\inf_{0\leq t\leq \delta}(\frac{S^n_{t+\tau}-n\widehat{s}_{t+\tau}}{a_n}-\frac{S^n_\tau-n\widehat{s}_\tau}{a_n})<-\epsilon,\notag\\ &\frac{n\delta_n }{a_n}\leq 1, \|\nu^n\|\leq M\big)=-\infty.
\end{align}
Equation \eqref{equ 5.7} follows from Equations \eqref{equ 5.9} and \eqref{equ 5.10} directly.

Let $e_2=(0,1)^T$ and $\chi_2(t)\equiv e_2$, then Equation \eqref{equ 5.8} follows from a similar analysis with that leading to Equation \eqref{equ 5.7} and the proof is complete.

\qed

\quad

\textbf{Acknowledgments.} The author is grateful to the financial
support from the National Natural Science Foundation of China with
grant number 11501542.

{}
\end{document}